\documentclass[11pt]{amsart}

\usepackage[english]{babel} 
\usepackage[T1]{fontenc}
\usepackage[utf8]{inputenc}

\usepackage{indentfirst}
\usepackage[left=2.5cm, right=2.5cm, top=2cm, bottom=2cm]{geometry} 
\usepackage{amsmath,amssymb,mathrsfs} 
\usepackage{graphicx}  
\usepackage{float, caption} 
\usepackage{wrapfig} 
\usepackage{framed}
\usepackage{array, multirow, tabularx} 
\usepackage{tikz} 
\usepackage{pifont} 
\usetikzlibrary{matrix,arrows,decorations.pathmorphing,patterns}
\usepackage{mathtools}

\usepackage[bookmarksnumbered=true,pdfstartview=FitH,pdftex,pdfborder={0 0 0},linkcolor=blue,citecolor=blue,colorlinks=true]{hyperref} 

\usepackage{multido}

{\newtheorem{de}{Definition}[section]}

{\newtheorem{theo}[de]{\textsc{Theorem}}}

{\newtheorem{prop}[de]{Proposition}}

{\newtheorem{lem}[de]{Lemma}}

{\newtheorem{cor}[de]{Corollary}}

{\theoremstyle{remark}
\newtheorem{remar}[de]{Remark}}

{\theoremstyle{remark}
}

{\theoremstyle{remark}
\newtheorem{ex}[de]{Example}}

{\theoremstyle{remark}
\newtheorem{nota}[de]{Notation}}

\newcommand{\etape}[1]{{\large\textit{#1}}}

\newcommand{\dr}{\ensuremath{\partial}}

\newcommand{\s}{\medskip}

\newcommand{\C}{\ensuremath{\mathbb{C}}}

\newcommand{\ssi}{if and only if}

\newcommand{\incl}{\ensuremath{\hookrightarrow}}
\newcommand{\ov}[1]{\ensuremath{\overline{#1}}}
\newcommand{\w}{\ensuremath{\omega}}
\newcommand{\co}{\ensuremath{\mathscr{O}}}
\newcommand{\wt}[1]{\ensuremath{\widetilde{#1}}}
\newcommand{\m}{\ensuremath{\mathfrak{m}}}
\newcommand{\f}[2]{\ensuremath{\frac{#1}{#2}}}

\newcommand{\mc}[1]{\ensuremath{\mathscr{#1}}}

\newcommand{\Cx}{\ensuremath{\C\{\underline{x}\}}}

\newcommand{\se}{\ensuremath{\subseteq}}

\newcommand{\lrp}[1]{\ensuremath{\left( #1\right)}}
\newcommand{\lra}[1]{\left\{#1\right\}}

\newcommand{\derl}{\ensuremath{\mathrm{Der}(-\log D)}}
\newcommand{\Cxy}{\ensuremath{\C\{x,y\}}}

\newcommand{\val}{\ensuremath{\mathrm{val}}}
\newcommand{\res}{\ensuremath{\mathrm{res}}}

\newcommand{\Z}{\ensuremath{\mathbb{Z}}}
\newcommand{\wo}{\ensuremath{w^{(0)}}}

\newcommand{\der}{\ensuremath{\mathrm{Der}}}
\newcommand{\undun}{\underline{1}}

\newcommand{\croix}[2]{\draw[thick] (#1+0.07,#2-0.07)-- (#1-0.07,#2+0.07);\draw[thick] (#1-0.07,#2-0.07)-- (#1+0.07,#2+0.07);}

\newcommand{\cercleg}[2]{\draw[thick,color=gray] (#1,#2) circle(0.13);}

\newcommand{\vide}{\emptyset}
\newcommand{\bs}{\backslash}
\newcommand{\dd}{\mathrm{d}}
\newcommand{\N}{\ensuremath{\mathbb{N}}}
\newcommand{\unp}{\left\{1,\ldots,p\right\}}
\newcommand{\unk}{\left\{1,\ldots,k\right\}}
\newcommand{\wh}[1]{\widehat{#1}}

\newcommand{\loga}{logarithmic}

\title{{On the values of logarithmic residues along curves}}
\author[D.~Pol]{Delphine Pol}
\address{
D.Pol\\
 Universit\'e d'Angers, D\'epartement de Math\'ematiques\\
 LAREMA, CNRS UMR n\textsuperscript{o}6093\\
 2 Bd Lavoisier\\
 49045 Angers\\
 France
 }
 \email{\href{pol@math.univ-angers.fr}{pol@math.univ-angers.fr}}
\date{\today}


\newcommand{\croixr}[2]{\draw[thick,color=gray] (#1+0.07,#2-0.07)-- (#1-0.07,#2+0.07);\draw[thick,color=gray] (#1-0.07,#2-0.07)-- (#1+0.07,#2+0.07);}

\def\restriction#1#2{\mathchoice
              {\setbox1\hbox{${\displaystyle #1}_{\scriptstyle #2}$}
              \restrictionaux{#1}{#2}}
              {\setbox1\hbox{${\textstyle #1}_{\scriptstyle #2}$}
              \restrictionaux{#1}{#2}}
              {\setbox1\hbox{${\scriptstyle #1}_{\scriptscriptstyle #2}$}
              \restrictionaux{#1}{#2}}
              {\setbox1\hbox{${\scriptscriptstyle #1}_{\scriptscriptstyle #2}$}
              \restrictionaux{#1}{#2}}}
\def\restrictionaux#1#2{{#1\,\smash{\vrule height .8\ht1 depth .85\dp1}}_{\,#2}}

\begin{document}

\begin{abstract}
We consider the germ of a reduced curve, possibly reducible. F.Delgado de la Mata proved that such a curve is Gorenstein if and only if its semigroup of values is symmetrical. We extend here this symmetry property to any fractional ideal of a Gorenstein curve. We then focus on the set of values of the module of \loga{} residues  along plane curves or complete intersection curves, which determines and is determined by the values of the Jacobian ideal thanks to our symmetry Theorem. Moreover, we give the relation with K\"ahler differentials, which are used in the analytic classification of plane branches. We also study the behaviour of \loga{} residues in an equisingular deformation of a plane curve.  
\end{abstract}

\maketitle

\section{Introduction}

 Let $D$ be the germ of a reduced hypersurface in $(\C^n,0)$ defined by $f\in\Cx:=\C\lra{x_1,\ldots,x_n}$ and with ring $\co_D=\Cx/(f)$. In his fundamental paper~\cite{saitolog}, K.Saito introduces the notions of logarithmic vector fields,  logarithmic differential forms and their residues. A logarithmic differential form is a meromorphic form on a neighbourhood of the origin in $\C^n$ which has simple poles along $D$  and such that its differential also has simple poles along $D$. A logarithmic $q$-form $\w$ satisfies: 
$$
g\w= \f{\dd f}{f}\wedge \xi+\eta,
$$
where $g\in\Cx$ does not induce a zero divisor in  $\co_D$, $\xi$ is a holomorphic $(q-1)$-form   and $\eta$ is a holomorphic $q$-form. Then, the logarithmic residue $\res^q(\w)$ of $\w$ is defined as the coefficient of $\frac{\dd f}{f}$, that is to say: $$\res^q(\w)=\restriction{\lrp{\frac{\xi}{g}}}{D}\in\Omega^{q-1}_D\otimes_{\co_D} Q(\co_D),$$ with $\Omega^{q-1}_D$ the module of K\"ahler differentials on $D$ and $Q(\co_D)$ the total ring of fractions of $\co_D$. We denote by $\mc{R}_D$ the $\co_D$-module of logarithmic residues of logarithmic $1$-forms. 

\s

In \cite{gsres}, M.Granger and M.Schulze prove that the $\co_D$-dual of the Jacobian ideal of $D$ is $\mc{R}_D$. If in addition $D$ is free, that is to say, if the module of logarithmic differential $1$-forms is a free $\Cx$-module, the converse also holds: the dual of $\mc{R}_D$ is the Jacobian ideal. They use this duality to prove a characterization of normal crossing divisors in terms of logarithmic residues: if the module $\mc{R}_D$ is equal to the module of weakly holomorphic functions on $D$ then $D$ is normal crossing in codimension~$1$. The converse implication was already proved in \cite{saitolog}.

\s

The purpose of this paper is to investigate more deeply the module of logarithmic residues. We focus on the case of plane curves or complete intersection curves. Plane curves are always free divisors, and they are the only singular free divisors with isolated singularities. The notion of multi-residues along complete intersections was introduced by A.G.Aleksandrov and A.Tsikh in \cite{alektsikh}.

\s

In order to describe the module of residues, or more generally, any fractional ideal, we will consider the set of values, which is defined as follows. Let $D=D_1\cup\dots\cup D_p$ be the germ of a reduced  complex analytic curve with $p$ irreducible components. The normalization of the local ring $\co_D$ induces a map $\val : Q(\co_D)\to \lrp{\Z\cup\lra{\infty}}^p$ called the value map, which associates with a fraction $g\in Q(\co_D)$ the $p$-uple of the valuations of $g$ along each irreducible component of $D$. Given a fractional ideal $I\subset Q(\co_D)$ (see Definition \ref{fracideal}), we denote by $\val(I)$ the set of values of the non zero divisors of $I$, and  $I^\vee:=\lra{g\in Q(\co_D); g\cdot I\subseteq \co_D}$ the dual of $I$.

\s

Let us explain the content of section \ref{part1}. We prove that the values of a fractional ideal and the values of its dual determine each other, and we give explicitly the relation between them. We then apply this result to the case of $\mc{R}_D$ and of the Jacobian ideal of $D$, denoted by $\mc{J}_D$, in parts \ref{part2} and~\ref{part3}. This relation is in fact a generalization of the following well-known theorem of Kunz in the case of irreducible curves:

\begin{theo}[\protect{\cite{kunzvalue}}]
\label{kunz}
If $D$ is irreducible, $\co_D$ is Gorenstein if and only if the  following property, which is called a symmetry property, is satisfied: for all $v\in\Z$, $$ v\in\val(\co_D)\iff \gamma-v-1\notin\val(\co_D),$$
where $\gamma$ is the conductor of $D$, that is to say, $\gamma=\min\lra{\alpha\in\N; \alpha+\N\subseteq \val(\co_D)}$.
\end{theo}
%We first show in section \ref{part1} that the values of the Jacobian ideal and the values of the module of \loga{} residues determine each other. The statement is easy for irreducible curves (see Remark~\ref{irred}):
%$$v\in\val(\mc{R}_D)\iff c-v-1\notin \val(\mc{J}_D)$$
%where $c\in\N$ is the conductor of the curve, i.e. $c=\min\lra{c\in\N; c+\N\subseteq \val(\co_D)}$.
%It is a generalization of the well-known symmetry of the semigroup for Gorenstein curves (see for example \cite{kunzvalue}).  
In \cite[Theorem 2.8]{delgado}, F.Delgado de la Mata generalizes the former result to the case of reducible curves. He proves that a curve is Gorenstein if and only if $\val(\co_D)$ satisfies a  symmetry property described below.

We prove here that Delgado's symmetry has an analogue which links the values of a fractional ideal and the values of its dual.   Whereas the symmetry is immediate for irreducible curves, the proof of this generalization of Delgado's Theorem is much more subtle. It  leads to the main result of this section, namely Theorem~\ref{symmetry-values}, which generalizes Theorem 2.4 of \cite{polcras} to any Gorenstein curve and any fractional ideal. To give the statement of our symmetry Theorem, we introduce the following notation for $i\in\unp$, $v\in\Z^p$ and $I$ a fractional ideal (see Notation~\ref{Delta}): $$\Delta_i(v,\val(I)):=\lra{\alpha\in\val(I); \alpha_i=v_i, \forall j\neq i, \alpha_j>v_j},$$ and $\Delta(v,\val(I))=\bigcup_{i=1}^p \Delta_i(v,\val(I))$. We consider the product order on $\Z^p$, that is to say, for $\alpha,\beta\in\Z^p$, $\alpha\leqslant\beta$ means that for all $i\in\lra{1,\ldots,p}, \alpha_i\leqslant \beta_i$. The conductor of $D$ is: $\gamma:=\inf\lra{\alpha\in\N^p; \alpha+\N^p\subseteq \val(\co_D)}$. We set $\underline{1}=(1,\ldots,1)$.   The statement of our main Theorem is: 
\begin{theo}
 \label{symmetry-values}
 Let $D$ be the germ of a reduced analytic curve with $p$ irreducible components.  Then, the ring $\co_D$ of the curve is a Gorenstein ring if and only if for all fractional ideals $I\subset Q(\co_D)$ the following property is satisfied for all $v\in\Z^p$:
 \begin{equation}
 \label{symmetry-values-eq}
 v\in \val(I^\vee)\iff \Delta(\gamma-v-\undun,\val(I))=\emptyset.
 \end{equation}
\end{theo} 
 Delgado's Theorem concerns the case $I=\co_D$ (see Theorem \ref{delgadotheo}). % has already proved in \cite{delgado} that the previous symmetry is satisfied by $\co_D$ if and only if the curve is Gorenstein. 
A similar symmetry was recently proved in \cite{kts} using  combinatorial methods which involve canonical modules.

 \s
 
 In part \ref{part2} and \ref{part3}, we use Theorem \ref{symmetry-values} to study the module of logarithmic residues along complete intersection curves, with a particular attention to the case of plane curves.

\s 
 
 In subsection~\ref{zero divisors}, we give some properties of the set of values of the module of \loga{} residues and of the Jacobian ideal for plane curves. We  investigate the zero divisors of $\mc{R}_D$ and $\mc{J}_D$, which are described in Propositions \ref{zerodiv:res} and \ref{zero:div:jac}. %They respectively come from the logarithmic residues and the logarithmic vector fields. %We first give a general property of fractional ideals:  the set of values of the non zero divisors determine the set of values of the zero divisors (see Proposition~\ref{zero:div}). We  prove that the zero divisors of $\mc{R}_D$ are determined by the logarithmic residues of the different possible unions of branches of $D$ (see Proposition \ref{zerodiv:res}), and the zero divisors for $\mc{J}_D$ comes from the logarithmic vector fields of the branches (see Proposition \ref{zero:div:jac}). 
 Thanks to our symmetry Theorem, we then determine  the conductor of $\mc{R}_D$, which is $-(m^{(1)},\ldots,m^{(p)})+\undun$, where $m^{(i)}$ is the multiplicity of the branch $D_i$. We  also mention the relation between logarithmic differential forms and the torsion of K\"ahler differentials. Thanks to this relation, we recover the result of O.Zariski on the dimension of the torsion of K\"ahler differentials, which is equal to the Tjurina number.

\s

 In subsection~\ref{complete-int}, we recall the theory of of multi-logarithmic differential forms and multi-residues along a reduced complete intersection  developed by  A.G.Aleksandrov and A.Tsikh in \cite{alektsikh}. Since our symmetry Theorem  is true for any Gorenstein curve, it is in particular true for  complete intersection curves. As in the hypersurface case, we again have a duality between multi-residues and the Jacobian ideal, so that their sets of values determine each other. Moreover, we prove here the following proposition (see Proposition~\ref{jc:kahler} for a more precise statement): %relation between the Jacobian ideal $\mc{J}_C$ and the K\"ahler differentials $\Omega^1_C$ of a complete intersection curve $C$: 
\begin{prop} 
Let $\Omega^1_C$ be the module of K\"ahler differentials along a reduced complete intersection curve $C$. The values of $\mc{J}_C$ and the values of $\Omega^1_C$ satisfy:
 $$\val(\mc{J}_C)=\gamma+\val(\Omega^1_C)-\undun.$$ 
 \end{prop}
 
 The set of values of K\"ahler differentials  is a major ingredient used in 
 \cite{hefez} and \cite{Hefez2} to study the problem of the analytic classification of plane curves with one or two branches. 
 
\s

 The last section is devoted to the study of the behaviour of \loga{} residues in an equisingular deformation of a plane curve. In particular, we define a stratification by the values of the \loga{} residues, which is the same as the stratification by K\"ahler differentials thanks to subsection~\ref{complete-int}. We prove that this stratification is finer than the stratification by the Tjurina number. We give an example in which the stratification by logarithmic residues is  strictly finer than the stratification by the Tjurina number (see Example~\ref{res:tau}). We show that the stratification by logarithmic residues is finite and constructible (see Propositions~\ref{finite} and \ref{construc}). We also give an example in which the frontier condition is not satisfied (see Example \ref{front}). 

\bigskip

\textbf{Acknowledgments.} The author is grateful to Michel Granger for many helpful discussions on the subject and his suggestion to use the result of Ragni Piene in the proof of Proposition~\ref{gamma}, and to Pedro Gonz\'alez-P\'erez and Patrick Popescu-Pampu for pointing out the papers of A.Hefez and M.E. Hernandes on the analytic classification of plane curves. The author also thanks Philipp Korell, Laura Tozzo and Patrick Popescu-Pampu for their careful reading of the previous version of this paper. The author thanks the anonymous referee for helpful comments and careful reading.

\section{The symmetry of values}
\label{part1}

This section is devoted to the main Theorem~\ref{symmetry-values}, which is a generalization of the symmetry Theorem 2.8 of \cite{delgado}. 

\s

We first recall several properties of fractional ideals and the notion of conductor of a curve. We then introduce some definitions and notations inspired by \cite{delgado} which appear in the statement of the main Theorem~\ref{symmetry-values}. We give in subsection \ref{proofth} a detailed proof of this Theorem, and then a property of the Poincar\'e series associated with a fractional ideal of a Gorenstein curve (see Proposition~\ref{sym:poly}).

\subsection{Properties of fractional ideals}

We  recall some basic results on the fractional ideals of a curve. In particular, we give several properties of the set of values of a fractional ideal, and we define its dual. %We also introduce the conductor of a curve  and some notations inspired by \cite{delgado}, which are essential to state our main Theorem \ref{symmetry-values}. 

\s

Let $D$ be the germ of  a reduced complex analytic curve, with $p$ irreducible components $D_1,\ldots,D_p$. The ring $\co_{D_i}$ of the branch $D_i$ is a one-dimensional integral domain, so that its normalization $\co_{\wt{D}_i}$ is isomorphic to $\C\lra{t_i}$ (see for example \cite[Corollary 4.4.10]{dejong}). By the splitting of normalization (see \cite[Theorem 1.5.20]{dejong}), the ring $\co_{\wt{D}}$ of the normalization of $D$ is $\co_{\wt{D}}=\bigoplus_{i=1}^p \C\lra{t_i}$. Moreover, the total rings of fractions $Q(\co_D)$ of $\co_D$ and $Q(\co_{\wt{D}})$ of $\co_{\wt{D}}$ are equal (see \cite[Exercise 4.4.16]{dejong}).  We then have : 
$$Q(\co_{\wt{D}})=Q(\co_D)=\bigoplus_{i=1}^p Q(\C\lra{t_i}).$$

\begin{de}
Let $g\in Q(\co_D)$. We define the valuation of $g$ along the branch $D_i$ as the order of $t_i$ of the image of $g$ by the surjection map $Q(\co_D)\to Q(\C\lra{t_i})$.  We denote the valuation of $g$ along $D_i$ by $\val_i(g)\in\Z\cup \lra{\infty}$, with the convention $\val_i(0)=\infty$.

We then define the value of $g$ by $\val(g):=\lrp{\val_1(g),\ldots,\val_p(g)}\in \lrp{\Z\cup\lra{\infty}}^p$.  
\end{de}
%Some authors write also $Q(\C\lra{t_i})=\C\lra{t_i}\left[\frac{1}{t_i}\right]$. 
\begin{de}
\label{fracideal}
Let $I\subset  Q(\co_D)$ be an $\co_D$-module. We call $I$ a fractional ideal if $I$ is of finite type over $\co_D$ and if $I$ contains a non zero divisor of $Q(\co_D)$. 
 We set:
$$\val(I):=\lra{\val(g); g\in I \text{ non zero divisor }}\subset \Z^p$$ and $$\ov{\val(I)}:=\lra{\val(g); g\in I}\subset \lrp{\Z\cup \lra{\infty}}^p.$$
\end{de}
\begin{remar}
We will prove in section~\ref{part2} that for a fractional ideal $I$, the set $\val(I)$ determines the set $\ov{\val(I)}$ (see Proposition~\ref{zero:div}).
\end{remar}

\begin{de}
Let $I\subset Q(\co_D)$ be a fractional ideal. The dual of $I$ is:
$$I^\vee:=\mathrm{Hom}_{\co_D}\lrp{I,\co_D}.$$
\end{de}

%\begin{de}
%Let $I\subset Q(\co_D)$ be a fractional ideal. The dual of $I$ is:
%$$I^\vee:=\lra{g\in Q(\co_D); g\cdot I\subseteq \co_D}$$
%\end{de}

\begin{remar}
We also have $I^\vee\simeq \lra{g\in Q(\co_D); g I\subseteq \co_D}$ (see for example \cite[Proof of Lemma 1.5.14]{dejong}).

%An element $g\in I^\vee$ induces a map $\varphi_g\in\mathrm{Hom}_{\co_D}(I,\co_D)$ defined as the multiplication by $g$, so that we have a natural map $\Psi : I^\vee\to\mathrm{Hom}_{\co_D}(I,\co_D)$. This map is in fact an isomorphism (see for example \cite[Proof of Lemma 1.5.14]{dejong}). 
\end{remar}

%From the definition of $I^\vee$, we see that 

\begin{lem}
\label{lem:dual}
The dual $I^\vee$ of a fractional ideal $I$ is also a fractional ideal. Moreover, if $I,J$ are two fractional ideals satisfying $J\subseteq I$, then $J^\vee\supseteq I^\vee$. 
\end{lem}

\s

\begin{de}
The conductor ideal of the curve $D$ is $\mc{C}_D:=\co_{\wt{D}}^\vee$.
\end{de}

For $\alpha\in\Z^p$, we set $t^\alpha:=(t_1^{\alpha_1},\ldots,t_p^{\alpha_p})\in Q(\co_D)$. The conductor ideal $\mc{C}_D$ is a fractional ideal of $\co_D$, and it is also an ideal in $\co_{\wt{D}}$. It implies the following property:
\begin{lem}
There exists $\gamma\in\N^p$ such that $\mc{C}_D=t^\gamma\co_{\wt{D}}$. This element $\gamma$ is called the conductor of the curve $D$. 
\end{lem}
%\begin{proof}The conductor ideal $\mc{C}_D$ is a fractional ideal of $\co_D$.  Moreover, it is also an ideal in $\co_{\wt{D}}=\bigoplus_{i=1}^p \C\lra{t_i}$, so that there exists $\gamma=(\gamma_1,\ldots, \gamma_p)\in\N^p$,  such that $\mc{C}_D=t^\gamma\co_{\wt{D}}$.
%\end{proof}

\s

We consider the  product order on $\mathbb{Z}^p$ we defined in the introduction. In particular, for $\alpha,\beta\in\mathbb{Z}^p$,   $\inf(\alpha,\beta)=\big(\min(\alpha_1,\beta_1),\ldots,\min(\alpha_p,\beta_p)\big)$.  %We set $\underline{1}=(1,\ldots,1)$.  
 
The conductor $\gamma$ satisfies: \begin{equation}
 \gamma=\inf\lra{\alpha\in\N^p; \alpha+\N^p\subseteq \val(\co_D)}.
 \end{equation}

We will need the following properties, which should be compared with  \cite[1.1.2, 1.1.3]{delgado}:

\begin{prop}
\label{propinf}

For a fractional ideal $I\subset Q(\co_D)$, if $v\in\val(I)$ and $v'\in\ov{\val(I)}$, then $\inf(v,v')\in\val(I)$. 

Similarly, if $v,v'\in\ov{\val(I)}$, then $\inf(v,v')\in\ov{\val(I)}$. 
\end{prop}
%\begin{remar}
%If $v\in\val(I)$ and $v'\in \ov{\val(I)}$, then $\inf(v,v')\in\val(I)$.
%\end{remar}

\begin{prop}
\label{valquimonte}
Let $v\neq v'\in\val(I)$. If there exists $i\in\unp$ such that $v_i=v'_i$, then there exists $v''\in\val(I)$ such that: 
\begin{enumerate}
\item $v''_i>v_i$
\item  For all $j\in\lra{1,\ldots,p}$, $v''_j\geqslant \min(v_j,v'_j)$
\item For all $j\in\lra{1,\ldots,p}$ such that $v_j\neq v'_j$, we have $v''_j=\min(v_j,v'_j)$
\end{enumerate}
\end{prop}

Proposition~\ref{propinf} is a consequence of the fact that the value of a general linear combination of two elements is equal to the minimum for the product order of the values of these two elements. Proposition~\ref{valquimonte} comes from the fact that a convenient linear combination will increase the valuation on the component $D_i$, but we cannot say what happens on the other components where the equality holds.

\s

From the definition of a fractional ideal, we have the following inclusions, which will be useful in subsection \ref{proofth}:
\begin{lem}
\label{incl}
Let $I$ be a fractional ideal. Then there exist $\nu$ and $\lambda$ in $\Z^p$ such that 
\begin{equation}
\label{incleq}
t^\nu \co_{\wt{D}}\subseteq I\subseteq t^\lambda \co_{\wt{D}}.
\end{equation}
\end{lem}
In particular, it implies that $\nu+\N^p\subseteq \val(I)\subseteq \lambda+\N^p$. Moreover, if $\lambda'\leqslant \lambda$ and $\nu'\geqslant \nu$, we can replace in \eqref{incleq} $\lambda$ by $\lambda'$ and $\nu$ by $\nu'$.

By dualizing \eqref{incleq}, we obtain, since $\co_{\wt{D}}^\vee=\mc{C}_D=t^\gamma\co_{\wt{D}}$:

\begin{equation}
\label{incl:dual}
t^{\gamma-\lambda}\co_{\wt{D}}\subseteq I^\vee\subseteq t^{\gamma-\nu}\co_{\wt{D}}.
\end{equation}

The following proposition is a  key ingredient for the proof of Theorem \ref{symmetry-values} in the Gorenstein case. 

\begin{prop}[\protect{\cite[Theorem 21.21]{eisenbudalgebra},  \cite[Lemma 5.2.8]{dejong}}]
\label{duality:MCM}
Let $I\subset Q(\co_D)$ be a fractional ideal of a Gorenstein curve. Then: 
\begin{itemize}
\item We have: $I^{\vee\vee}=I$.
\item If $I\subset J$ are fractional ideals, $\dim_\C J/I=\dim_\C I^\vee/J^\vee$.
\end{itemize}
\end{prop}

 \s

%\subsection{Statement of the symmetry Theorem}

We end this subsection with the following notations, which are analogous to the notations of \cite{delgado}. 

\begin{nota}
\label{Delta}
Let us consider an arbitrary subset $\mc{E}$ of $\mathbb{Z}^p$ and $v\in\mathbb{Z}^p$. For $i\in\{1,\ldots,p\}$, we define:
$$\Delta_i(v,\mc{E})=\{\alpha\in\mc{E};\  \alpha_i=v_i \text{ and } \forall j\neq i, \alpha_j>v_j\},$$
and $\Delta(v,\mc{E})=\bigcup_{i=1}^p \Delta_i(v,\mc{E})$. For a fractional ideal $I\subset Q(\co_D)$, we write $\Delta(v,I)$ instead of $\Delta(v,\val(I))$.
\end{nota}

The following picture illustrates the case $p=2$. Let us consider the subset $\mc{E}$ of $\Z^2$ defined by all the crosses. The grey crosses correspond to the elements of $\Delta(v,\mc{E})$, that is to say, $\Delta(v,\mc{E})=\lra{(3,1),(4,1),(2,3),(2,4),(2,5)}$.% $\mc{E}=\lra{(1,0)}\cup\lra{((i,1), (i,3),i\in\lra{2,3,4}}\cup \lra{(2,4),(2,5),(3,4)}$ and $v=(2,1)$.The crosses represent the elements of~$\mc{E}$. The grey ones are the elements of $\Delta(v,\mc{E})$.

\begin{figure}[H]
\begin{center}
\begin{tikzpicture}
%\shorthandoff{:} 
\draw [step=0.5, gray, very thin] (0,0.5) grid (3,3);
\draw [->, thick] (0,0.5) -- (3.5,0.5) node[at end, below] {{\small $\val_1$}};
\draw [->, thick] (0,0.5) -- (0,3.5) node[at end, left] {{\small $\val_2$}};
\croix{1}{1};
\croixr{1}{2};
\croixr{1.5}{1};
\croixr{2}{1};
\croix{1.5}{2};
\croix{2}{2};
\croix{1.5}{2.5};
\croix{0.5}{0.5};
\croixr{1}{2.5};
\croixr{1}{3};
\draw (1,1) node[below]{$v$};
\draw (0,0.5) node[below]{$0$};
\draw (2.5,0.5) node[below]{$5$};
\draw (0,3) node[left]{$5$};
%\draw[-,thick,red] (1,1.5)--(1,3.2);
%\draw[-,thick,red] (1.5,1)--(3.2,1); 
%\draw[->] (7,1)--(10,1);
%\draw[->] (7,1)--(7,3.5);
%\draw[->] (7,1)--(6,-0.5);
%\draw (7,1) node[below]{$v$};
%\draw[-,gray] (7.5,1.5)--(7.5,3.5);
%\draw[-,gray] (7.5,1.5)--(10,1.5);
%\draw[-,gray] (7,0.5)--(6,-1);
%\draw[-,gray] (7,0.5)--(10,0.5);
\end{tikzpicture}
\end{center}

\caption{$\Delta(v,\mc{E})$ for $p=2$}
\end{figure}

We recall here the statement of Delgado's Theorem:

\begin{theo}[\protect{\cite[Theorem 2.8]{delgado}}]
\label{delgadotheo}
 Let $D$ be the germ of a reduced curve with $p$ irreducible components.  Then, the ring $\co_D$ of the curve is a Gorenstein ring if and only if for all $v\in\Z^p$,
 \begin{equation}
 v\in \val(\co_D)\iff \Delta(\gamma-v-\undun,\co_D)=\emptyset.
 \end{equation}
\end{theo}

\subsection{Proof of Theorem~\ref{symmetry-values}}
\label{proofth}

The proof of Theorem \ref{symmetry-values} is developed in several steps. We first mention the implications which are easy consequences of \cite{delgado} (see Lemma \ref{goren} and Proposition \ref{sens1}). The remain of this subsection is then devoted to the missing implication, which needs much more work (see subsections \ref{steps} and \ref{fin}). 

\s

The following Lemma is a direct consequence of Theorem \ref{delgadotheo}.  Indeed, if the condition \eqref{symmetry-values-eq} is satisfied for all fractional ideals $I$, it is in particular satisfied by $\co_D$.
\begin{lem}
\label{goren}
Let $D$ be a reduced curve. 
If \eqref{symmetry-values-eq} is satisfied for all fractional ideals $I\subset Q(\co_D)$, then $\co_D$ is Gorenstein.
\end{lem}

\begin{remar}
It is not sufficient to check if \eqref{symmetry-values-eq} is satisfied for one fractional ideal $I$ to prove that the curve is Gorenstein. Indeed, by definition, for every curve, the equivalence~\eqref{symmetry-values-eq} is satisfied by $I=\co_{\wt{D}}$ and $I^\vee=\mc{C}_D$.  
\end{remar}

Our purpose now is to prove that for a Gorenstein curve and a fractional ideal $I\subset Q(\co_D)$, property \eqref{symmetry-values-eq} is satisfied. Nevertheless, some of the properties we will prove or mentioned are also satisfied by non Gorenstein curves, so that we first consider an arbitrary reduced curve $D$.

Let $I\subset Q(\co_D)$ be a fractional ideal. To prove the implication $\Rightarrow$ of \eqref{symmetry-values-eq}, we need the following result: 

\begin{prop}[\protect{\cite[Corollary 1.9]{delgado}}]
\label{cor1.9del}
Let $D$ be a reduced curve. Then: $$\Delta(\gamma-\undun,\co_D)=\emptyset.$$
\end{prop}

\begin{prop}
\label{sens1}
Let $D$ be a reduced curve. Let $v\in\Z^p$ and $I\subset Q(\co_D)$ be a fractional ideal. Then:
$$v\in \val(I^\vee)\Rightarrow \Delta(\gamma-v-\undun,I)=\emptyset.$$
\end{prop}
\begin{proof}
Let $v=(v_1,\ldots,v_p)\in \val(I^\vee)$, and $g\in I^\vee$ with $v=\val(g)$. We assume $\Delta(\gamma-v-\undun,I)\neq \emptyset$. For the sake of simplicity, we may assume that $\Delta_1(\gamma-v-\undun,I)\neq \emptyset$. It means that there exists $h\in I$ with $\val(h)=(\gamma_1-v_1-1,w_2,\ldots,w_p)\in \val(I)$ such that for all $j\geqslant 2$, $w_j> \gamma_j-v_j-1$. Since $gh\in\co_D$, we have $(\gamma_1-1,w_2+v_2,\ldots,w_p+v_p)\in\val(\co_D)$, with $w_j+v_j\geqslant \gamma_j$. Therefore, $\Delta_1(\gamma-\undun,\co_D)\neq\emptyset$. 

Nevertheless, by Proposition \ref{cor1.9del}, $\Delta(\gamma-\undun,\co_D)=\emptyset$, which leads to a contradiction. Therefore, $\Delta(\gamma-v-\undun,I)=\emptyset$.
\end{proof}

\begin{nota}
\label{V}
We set $\mc{V}=\lra{v\in\Z^p; \Delta(\gamma-v-\undun,I)=\emptyset}$. 
\end{nota}

Proposition~\ref{sens1} tells us that the set $\mc{V}$ contains the values of $I^\vee$, but \textit{a priori} it may be bigger. In particular, it is not obvious that $\mc{V}$ is the set of values of an $\co_D$-module. Our purpose here is to prove that $\mc{V}$ is indeed equal to $\val(I^\vee)$ when $D$ is Gorenstein.

\begin{prop}
\label{irred}
Let $D$ be an irreducible Gorenstein curve. The statement of Theorem~\ref{symmetry-values} can be rephrased as follows: for all $v\in\Z$, $v\in\val(I^\vee)$ if and only if $\gamma-v-1\notin\val(I)$.
\end{prop}
\begin{proof}
This proposition is a generalization of Kunz's theorem~\ref{kunz}.  We have: \begin{align*}
\dim_\C I^\vee/t^{\gamma-\lambda}\C\lra{t}&=\mathrm{Card}\lrp{\val(I^\vee)\cap \lrp{\gamma-\lambda+\N}^c}\\
\dim_\C t^\lambda\C\lra{t}/I&=\mathrm{Card}\big((\lambda+\N)\cap \lrp{\val(I)}^c\big)
\end{align*} Since by Proposition~\ref{duality:MCM}, $\dim_\C I^\nu/t^{\gamma-\lambda}\C\lra{t}=\dim_\C t^\lambda\C\lra{t}/I$,  we have the result. 
\end{proof}

\medskip

The proof for a reducible Gorenstein curve is based on a more intricate dimension argument.

\subsubsection{Dimension and values}

As we have seen in Proposition \ref{irred}, it is easy to compute dimensions from the set of values in the irreducible case. The purpose of this subsection is to give a combinatorial method to compute some dimensions from the set of values. 

\s

Let $v\in\Z^p$. We set $I_v=\lra{g\in I; \val(g)\geqslant v}$ and  $\ell(v,I)=\dim_\C I/I_v$. Since $\co_D$ is one-dimensional, we have $\ell(v,I)<\infty$. 

We denote by $(e_1,\ldots,e_p)$ the canonical basis of $\Z^p$. For $\mc{E}\subseteq \mathbb{Z}^p$ and $v\in\mathbb{Z}^p$, let 
$$\Lambda_i(v,\mc{E})=\{\alpha\in\mc{E}\ ;\ \alpha_i=v_i \text{ and } \alpha\geqslant v\}.$$

We have the following inclusion: $\Delta_i(v,\mc{E})\subseteq \Lambda_i(v,\mc{E})$. 
We then have  (see \cite[Proposition 1.11]{delgado}):

\begin{prop}
\label{dimension}
For all $v\in \mathbb{Z}^p$, $\ell(v+e_i,I)-\ell(v,I)=\dim_\C I_v/I_{v+e_i}\in\{0,1\}$ and moreover $\ell(v+e_i,I)=\ell(v,I)+1$ if and only if $\Lambda_i(v, \val(I))\neq\emptyset$.
\end{prop}

Thanks to this proposition we can compute some dimensions from the set of values: 

\begin{cor}
\label{suite}
Let $\nu,\lambda\in\Z^p$ such that $\nu+\N^p\subseteq \val(I)\subseteq \lambda+\N^p$ (see Lemma~\ref{incl}). Let $(\alpha^{(j)})_{0\leqslant j\leqslant M+1}$ be a finite sequence of elements of $\Z^p$  with $M=-1+\sum_{i=1}^p (\nu_i-\lambda_i)$,  which satisfies:
\begin{itemize}
\item $\alpha^{(0)}=\lambda$ and $\alpha^{(M+1)}=\nu$
\item For all $j\in\{0,\ldots,M\}$, there exists $i(j)\in\{1,\ldots,p\}$ such that $\alpha^{(j+1)}=\alpha^{(j)}+e_{i(j)}$
\end{itemize}
Then: 
\begin{equation}
\label{dim:seq}
\dim_\C I/t^\nu\co_{\wt{D}}=\ell(\nu,I)=\mathrm{Card}\lra{j\in\{0,\ldots,M\} ; \Lambda_{i(j)}(\alpha^{(j)},\val(I))\neq\emptyset}.
\end{equation}
\end{cor}

\begin{ex}
\label{example:corollary}
The following example illustrates Corollary~\ref{suite} for $p=2$. We consider the plane curve $D$ defined by $f(x,y)=(x^2-y^3)(x^4-y^3)$. A parametrization of this place curve is given by $x=(t_1^3,t_2^3)$, $y=(t_1^2,t_2^4)$. We consider the Jacobian ideal $I=\mc{J}_D$ of $D$.  In particular, $\val\lrp{\frac{\dr f}{\dr x}}=(9,15)$ and $\val\lrp{\frac{\dr f}{\dr y}}=(10,14)$. Therefore, by Proposition~\ref{propinf}, we have $(9,14)\in\val(\mc{J}_D)$.  We represent by crosses the elements of $I$. 

We can choose for example $\lambda=(8,13)$ and $\nu=(13,21)$. We consider the sequence $\alpha$ defined by the grey circles on figure~\ref{val:jd:morpi}. In particular, $\alpha^{(0)}=\lambda$ and $\alpha^{(13)}=\nu$. 
The sets $\Lambda_{i(j)}(\alpha^{(j)},\val(I))$ for $j\in\lra{0,\ldots,12}$ corresponds to the crosses which are on the thick black lines. By the corollary, for this example, we have: $\dim_\C I/t^\nu\co_{\wt{D}}=7$.

\begin{figure}[H]
\begin{center}
\begin{tikzpicture}
%\shorthandoff{:} 
\draw [step=0.5, gray, very thin] (0,0.5) grid (4.5,5.5);
\draw [->] (0,0.5) -- (5,0.5) node[at end, below] {{\small $\val_1$}};
\draw [->] (0,0.5) -- (0,6) node[at end, left] {{\small $\val_2$}};

\croix{1}{1.5}
\croix{1}{2}
\croix{1.5}{1.5}
\croix{2}{4}
\multido{\roo=2+0.5}{6}{\croix{\roo}{3}}
\multido{\roo=2+0.5}{6}{\croix{\roo}{3.5}}
\multido{\roo=2.5+0.5}{5}{\croix{\roo}{4.5}}
\multido{\roo=2.5+0.5}{5}{\croix{\roo}{5}}
\multido{\roo=2.5+0.5}{5}{\croix{\roo}{5.5}}
\cercleg{0.5}{1}
\cercleg{1}{1}
\cercleg{1.5}{1}
\cercleg{1.5}{1.5}
\cercleg{1.5}{2}
\cercleg{2}{2}
\cercleg{2}{2.5}
\cercleg{2}{3}
\cercleg{2}{3.5}
\cercleg{2.5}{3.5}
\cercleg{2.5}{3.5}
\cercleg{3}{3.5}
\cercleg{3}{4}
\cercleg{3}{4.5}
\cercleg{3}{5}
\draw (3,5.2) node[right]{$\nu$};
\draw (0.5,1) node[below]{$\lambda$};
\draw[-,very thick] (0.5,1)--(0.5,5.5);
\draw[-,very thick] (1,1)--(1,5.5);
\draw[-,very thick] (1.5,1)--(4.5,1);
\draw[-,very thick] (1.5,1.5)--(4.5,1.5);
\draw[-,very thick] (1.5,2)--(1.5,5.5);
\draw[-,very thick] (2,2)--(4.5,2);
\draw[-,very thick] (2,2.5)--(4.5,2.5);
\draw[-,very thick] (2,3)--(4.5,3);
\draw[-,very thick] (2,3.5)--(2,5.5);
\draw[-,very thick] (2.5,3.5)--(2.5,5.5);
\draw[-,very thick] (3,3.5)--(4.5,3.5);
\draw[-,very thick] (3,4)--(4.5,4);
\draw[-,very thick] (3,4.5)--(4.5,4.5);
%{\multido{\roo=0+0.5}{3}{\cercler{0}{\roo}}
%\multido{\roo=0+0.5}{3}{\cercler{0.5}{\roo}}
%\cercler{1}{2}
%\cercler{2}{4}
%\multido{\roo=0+0.5}{4}{\cercler{\roo}{2.5}}
%\multido{\roo=1.5+0.5}{3}{\cercler{1.5}{\roo}}}

\draw (1.5,0.5) node[below]{$10$};
\draw (4,0.5) node[below]{$15$};
\draw (-0.05,2) node[left]{$15$};
\draw (-0.05,4.5) node[left]{$20$};
%\draw (0.5,0.5) node[right]{$\gamma$};
%\cerclef{0.5}{0.5}
\end{tikzpicture}
\end{center}
\caption{Illustration of Corollary~\ref{suite}}
\label{val:jd:morpi}
\end{figure}
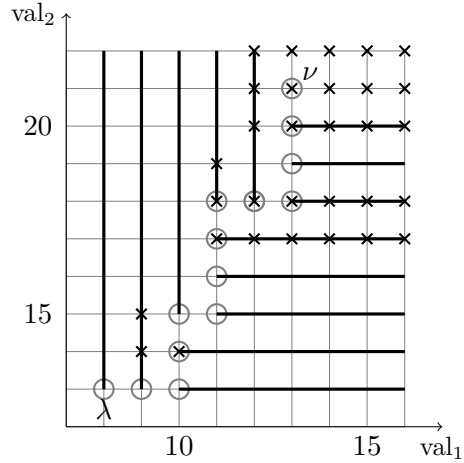
\end{ex}

\subsubsection{Preliminary steps}
\label{steps}

We recall that $\mc{V}=\lra{v\in \Z^p; \Delta(\gamma-v-\undun,I)=\emptyset}$, and this set contains $\val(I^\vee)$. The purpose of this section is to show that if the inclusion $\val(I^\vee)\subseteq \mc{V}$ is strict, then it has some combinatorial and numerical consequences (see Lemma \ref{lem:lp} and Proposition \ref{prop:cle}). We then prove in subsection \ref{fin} that these criteria lead to a contradiction in the Gorenstein case, which finishes the proof of Theorem~\ref{symmetry-values}. 

\medskip

\textit{First step}

We first show that if $\mc{V}\neq \val(I^\vee)$, then there is an element $w\in\mc{V}\bs \val(I^\vee)$ which satisfies some properties which will be used in the next steps. 

Let us assume that $\mc{V}\neq \val(I^\vee)$, and let $\wo\in\mc{V}\backslash \val(I^\vee)$ be "an intruder". By Lemma~\ref{incl}, there exist $\lambda,\nu\in\Z^p$ satisfying $\nu+\N^p\subseteq \val(I)\subseteq \lambda+\N^p$ and $\gamma-\nu\leqslant  \wo\leqslant \gamma-\lambda$.  For the remainder of the proof, we fix such $\lambda,\nu$.

The following proposition gives an essential property of  $\wo$:
\begin{prop}
\label{intrus}
There exists $j\in\{1,\ldots,p\}$ such that $\Lambda_j(w^{(0)},\val(I^\vee))=\vide$. Moreover, the corresponding coordinate satisfies $w^{(0)}_j<\gamma_j-\lambda_j$. 
\end{prop}

\begin{proof}
If for all $i\in\{1,\ldots,p\}$, $\Lambda_i(\wo,\val(I^\vee))\neq\vide$, then for all $i\in\{1,\ldots, p\}$ there exists $\alpha^{(i)}\in\val(I^\vee)$ such that $\alpha^{(i)}_i=\wo_i$ and $\alpha^{(i)}_j\geqslant \wo_j$. As a consequence, by Proposition~\ref{propinf}, $\inf(\alpha^{(1)},\ldots,\alpha^{(p)})=\wo\in\val(I^\vee)$, which is a contradiction. It gives the existence of a $j\in\{1,\ldots,p\}$ such that $\Lambda_j(w^{(0)},\val(I^\vee))=\vide$. It is immediate to see that $\wo_j<\gamma_j-\lambda_j$ since if $\wo_j=\gamma_j-\lambda_j$, then $\gamma-\lambda\in\Lambda_j(w^{(0)},\val(I^\vee))$, which contradicts the emptiness. 
\end{proof}

\medskip

\textit{Second step}

For the sake of simplicity, we assume that $\Lambda_p(\wo, \val(I^\vee))=\emptyset$.  Corollary~\ref{suite} together with a convenient finite sequence $\alpha$ can be used to compute the dimension of the quotient $I^\vee/t^{\gamma-\lambda}\co_{\wt{D}}$. We compare it with the number $\ell=\mathrm{Card}\lra{j\in\{0,\ldots,n_0-1\}\ ;\ \Lambda_{i(j)}(\alpha^{(j)},\mc{V})\neq\emptyset}$, which may \textit{a priori} depend on the chosen sequence $\alpha$.

\medskip

In order to compute $\dim_\C I^\vee/t^{\gamma-\lambda}\co_{\wt{D}}$, we consider a sequence $(\alpha^{(j)})_{0\leqslant j\leqslant n_0}$ with $n_0=\sum_{i=1}^p \lrp{\nu_i-\lambda_i}$ satisfying: 
\begin{itemize}
\item $\alpha^{(0)}=\gamma-\nu$ and $\alpha^{(n_0)}=\gamma-\lambda$
\item for all $j\in\lra{0,\ldots n_0-1}$, there exists $i(j)\in\lra{1,\ldots,p}$ such that $\alpha^{(j+1)}=\alpha^{(j)}+e_{i(j)}$
\item there exists $j_0\in \lra{0,\ldots,n_0-1}$ such that $\alpha^{(j_0)}=\wo$ and $\alpha^{(j_0+1)}=\wo+e_p$
\end{itemize}

The existence of such a sequence follows from Proposition~\ref{intrus}.  Moreover, this sequence satisfies the required properties of Corollary~\ref{suite}.

\medskip

Let us consider again the plane curve defined by $f(x,y)=(x^2-y^3)(x^4-y^3)$ and the ideal $I=\mc{J}_D$ of example~\ref{example:corollary}. By computing $\val(\co_D)$, one can see that the conductor $\gamma$ satisfies $\val(\gamma)=(8,12)$. The black crosses on figure~\ref{suite:alpha} represent the elements of $\mc{V}=\lra{w\in\Z^2, \Delta(\gamma-v-\undun,I)=\emptyset}$. Let us  assume for example that $\wo=(-2,-4)\notin\val(I^\vee)$. Then we can for instance consider the sequence $\alpha$ defined by the grey circles, where $\alpha^{(0)}=\gamma-\nu=(-5,-9)$ and $\alpha^{(n_0)}=\alpha^{(13)}=\gamma-\lambda=(0,-1)$.  %For example, when $p=2$, we can choose a sequence $\alpha$ as follows, where the elements of $\alpha$ are represented by grey circles and the elements of $\mc{V}$ by black crosses:

\begin{figure}[H]
\begin{center}
\begin{tikzpicture}
%\shorthandoff{:} 

\draw [step=0.5, gray, very thin] (-3.5,-5) grid (0.5,0.5);
\draw [->] (-4,0) -- (1,0) node[at end, below] {{\small $\val_1$}};
\draw [->] (0,-5.5) -- (0,1) node[at end, left] {{\small $\val_2$}};

\croix{-2}{-4}
\multido{\roo=-1.5+0.5}{5}{\croix{\roo}{-2.5}}
\croix{-1}{-2}
\multido{\roo=-2+0.5}{2}{\croix{-1.5}{\roo}}
\multido{\roo=-1+0.5}{4}{\croix{-0.5}{\roo}}
\multido{\roo=-1+0.5}{4}{\croix{0}{\roo}}
\multido{\roo=-1+0.5}{4}{\croix{0.5}{\roo}}
\draw (-0.6,-2) [above] node{$w^{(0)}$};

\cercleg{0}{-0.5}
\cercleg{-0.5}{-0.5}
\cercleg{-1}{-0.5}
\cercleg{-1}{-1}
\cercleg{-1}{-1.5}
\cercleg{-1}{-2}
\multido{\roo=-4.5+0.5}{5}{\cercleg{-1}{\roo}}
\cercleg{-1.5}{-4.5}
\cercleg{-2}{-4.5}
\cercleg{-2.5}{-4.5}
%\draw (-2.7,-4.5)[below] node{$\gamma-\nu$};
%\draw (0.2,-0.55) [above] node{$\gamma-\lambda$};
\draw (-2.5,0.05) node[above]{$-5$};
\draw (0.5,-2.5) node[right]{$-5$};
\draw (0.5,-5) node[right]{$-10$};
\draw (0.25,0.25) node{$0$};
\end{tikzpicture}
\end{center}
%\begin{center}
%\begin{tikzpicture}
%
%%\draw [step=0.5, gray, very thin] (-4.5,-3.5) grid (0.25,0.25);
%%\draw [->, thick] (-5,0) -- (1,0) node[at end, right] {$\val_1$};
%%\draw [->, thick] (0,-4) -- (0,1) node[at end, left] {$\val_2$};
%\draw [step=0.5, gray, very thin] (-4,-3.5) grid (0.25,0.25);
%\draw [->, thick] (-4.5,0) -- (1,0) node[at end, right] {$\val_1$};
%\draw [->, thick] (0,-4) -- (0,1) node[at end, left] {$\val_2$};
%
%\draw[fill=black] (-2,-1.5) circle(0.05) node[above]{$\wo$};
%%\croix{-4}{-3} 
%%\croix{0}{0}
%\draw (-3.5,-3) node[below]{$\alpha^{(0)}$};
%\draw (-3.5,-3) node[above]{$\gamma-\nu$};
%\multido{\roo=0+0.5}{4}{\cercleg{-3.5+\roo}{-3}}
%\multido{\roo=0+0.5}{7}{\cercleg{-2}{-3+\roo}}
%\multido{\roo=0+0.5}{5}{\cercleg{-2+\roo}{0}}
%\draw (0.5,-0.2) node{$\alpha^{(n_0)}$};
%\draw (0,0.3) node[right]{$\gamma-\lambda$};
%
%\end{tikzpicture}
%\end{center}
\caption{}
\label{suite:alpha}
\end{figure}

\begin{remar}
In \cite{polcras} we choose a sort of "canonical" sequence, but it is in fact unnecessary and the  presentation here is simpler. 
\end{remar}

From Corollary~\ref{suite},  we have: 
\begin{equation}
\label{dimIvC}
\dim_\C I^\vee/t^{\gamma-\lambda}\co_{\wt{D}}=\mathrm{Card}\lra{j\in\{0,\ldots,n_0-1\}\ ;\ \Lambda_{i(j)}(\alpha^{(j)},\val(I^\vee))\neq\emptyset}.
\end{equation}

We want to compare this dimension with the following number $\ell$:
\begin{equation}
\label{defellprime}
\ell=\mathrm{Card}\lra{j\in\{0,\ldots,n_0-1\}\ ;\ \Lambda_{i(j)}(\alpha^{(j)},\mc{V})\neq\emptyset}.
\end{equation}

\begin{lem}
\label{lem:lp}
For the sequence $\alpha$ defined above, we have:
\begin{equation}
\label{ellprime}
\ell\geqslant 1+\dim_\C I^\vee/t^{\gamma-\lambda}\co_{\wt{D}} .
\end{equation}
\end{lem}

\begin{proof}
 It is clear that $\Lambda_{i(j)}(\alpha^{(j)},\val(I^\vee))\neq \vide\Rightarrow \Lambda_{i(j)}(\alpha^{(j)},\mc{V})\neq \vide$. Moreover, since there exists $j_0$ such that $\alpha^{(j_0)}=\wo$ and $\alpha^{(j_0+1)}=\alpha^{(j_0)}+e_p$, $\Lambda_p(\alpha^{(j_0)},\mc{V})\neq\vide$, but the assumptions on  $\wo$ implies $\Lambda_p(\alpha^{(j_0)},\val(I^\vee))=\vide$. Hence the inequality.
\end{proof}

From now on, our sequence $\alpha$ is fixed. 

\medskip

\textit{Third step}

The purpose of this third step is to compare this number $\ell$ to $\dim_\C I/t^{\nu}\co_{\wt{D}}$.

 For $i\in\{0,\ldots,n_0\}$ we set $\beta^{(i)}=\gamma-\alpha^{(n_0-i)}$. The sequence  $\beta$ satisfies the properties of Corollary~\ref{suite}, so that it can be used to compute the dimension $\dim_\C I/t^\nu\co_{\wt{D}}$.

\medskip

We continue with the same example as before. For the sequence $\alpha$ of Figure \ref{suite:alpha}, we  represent the corresponding sequence $\beta$ on the following diagram. The sequence $\beta$ is defined for $i\in\lra{0,\ldots,13}$ by $\beta^{(i)}=\gamma-\alpha^{(13-i)}$. In particular, $\beta^{(0)}=\lambda$ and $\beta^{(13)}=\nu$. The elements of $\beta$ are represented by grey circles and the elements of $I$ by black crosses. In particular, $\beta^{(0)}=\lambda=(8,13)$ and $\beta^{(13)}=\nu=(13,21)$.  %, where $\vo=\gamma-\wo-\undun$:
\begin{figure}[H]
\begin{center}
\begin{tikzpicture}
%\shorthandoff{:} 
\draw [step=0.5, gray, very thin] (0,0) grid (4.5,5.5);
\draw [->] (0,0) -- (5,0) node[at end, below] {{\small $\val_1$}};
\draw [->] (0,0) -- (0,6) node[at end, left] {{\small $\val_2$}};

\croix{1}{1.5}
\croix{1}{2}
\croix{1.5}{1.5}
\croix{2}{4}
\multido{\roo=2+0.5}{6}{\croix{\roo}{3}}
\multido{\roo=2+0.5}{6}{\croix{\roo}{3.5}}
\multido{\roo=2.5+0.5}{5}{\croix{\roo}{4.5}}
\multido{\roo=2.5+0.5}{5}{\croix{\roo}{5}}
\multido{\roo=2.5+0.5}{5}{\croix{\roo}{5.5}}

\cercleg{0.5}{1}
\cercleg{1}{1}
\cercleg{1.5}{1}
\multido{\roo=1.5+0.5}{8}{\cercleg{1.5}{\roo}}
\multido{\roo=2+0.5}{3}{\cercleg{\roo}{5}}
\draw (0.5,1) [below] node{$\lambda$};
\draw (3.2,5)[above] node{$\nu$};
%{\multido{\roo=0+0.5}{3}{\cercler{0}{\roo}}
%\multido{\roo=0+0.5}{3}{\cercler{0.5}{\roo}}
%\cercler{1}{2}
%\cercler{2}{4}
%\multido{\roo=0+0.5}{4}{\cercler{\roo}{2.5}}
%\multido{\roo=1.5+0.5}{3}{\cercler{1.5}{\roo}}}

\draw (1.5,-0.05) node[below]{$10$};
\draw (4,-0.05) node[below]{$15$};
\draw (-0.05,2) node[left]{$15$};
\draw (-0.05,4.5) node[left]{$20$};
%\draw (0.5,0.5) node[right]{$\gamma$};
%\cerclef{0.5}{0.5}
\end{tikzpicture}
\end{center}
%\caption{Values of $\mc{J}_C$}
%\label{val:jd:morpi}
%\end{figure}
%\begin{center}
%\begin{tikzpicture}
%
%
%\draw [step=0.5, gray, very thin] (0,0) grid (4.5,4);
%\draw [->, thick] (0,0) -- (5,0) node[at end, below] {$\val_1$};
%\draw [->, thick] (0,0) -- (0,4.5) node[at end, left] {$\val_2$};
%
%\draw[fill=black] (2,1.5) circle(0.05) node[above]{$\vo$};
%%\croix{0.5}{0.5} 
%\draw (0.5,0.5) node[below]{$\beta^{(0)}$};
%\draw (0.5,0.5) node[above]{$\lambda$};
%%\croix{4.5}{3.5}
%\draw (4,3.4) node[below]{$\nu$};
%\multido{\roo=0+0.5}{5}{\cercleg{0.5+\roo}{0.5}}
%\multido{\roo=0+0.5}{7}{\cercleg{2.5}{0.5+\roo}}
%\multido{\roo=0+0.5}{4}{\cercleg{2.5+\roo}{3.5}}
%\draw (4,3.5) node[above]{$\beta^{(n_0)}$};
%
%\end{tikzpicture}
%\end{center}
\caption{}
\label{beta}
\end{figure}

The following proposition gives a relation between $\ell$ and $\dim_\C I/t^\nu\co_{\wt{D}}$:

\begin{prop}
\label{prop:cle}
With the above notations we have:
$$\ell\leqslant\sum_{i=1}^p{(\nu_i-\lambda_i)} -\dim_\C I/t^\nu\co_{\wt{D}}.$$
\end{prop}
To prove this proposition, we need the following lemma:

\begin{lem}
\label{lem:lam:v}
Let $w\in\Z^p$ and $i\in\unp$. Then:
$$\Lambda_i(w,\mc{V})\neq \emptyset \Rightarrow \Lambda_i(\gamma-w-e_i,\val(I))=\emptyset.$$
\end{lem}
\begin{proof}
Let $w'\in\Lambda_i(w,\mc{V})$. By the definition of $\mc{V}$, we have $\Delta(\gamma-w'-\undun,\val(I))=\emptyset$. Moreover, $(\gamma-w'-e_i)_i=\gamma_i-w_i-1$ and for $j\neq i$, $(\gamma-w'-e_i)_j=\gamma_j-w'_j\leqslant \gamma_j-w_j$.  Thus $\Lambda_i(\gamma-w-e_i,\val(I))=\Delta_i(\gamma-w'-\undun,\val(I))=\emptyset$, since $w'\in\mc{V}$.
\end{proof}

\textbf{Proof of Proposition~\ref{prop:cle}.}
We first notice that the two sequences $\alpha$ and $\beta$ have the same number of terms, namely $n_0+1=\sum_{i=1}^p (\nu_i-\lambda_i)\ +1$. 

By Corollary~\ref{suite}, we have:
\begin{equation}
\label{eq:prop:cle}
\dim_\C I/t^\nu\co_{\wt{D}}=\text{Card}\lra{j\in\{0,\ldots,n_0-1\}\ ;\  \Lambda_{i(n_0-j-1)}(\beta^{(j)},\val(I))\neq\vide}.
\end{equation}

We notice that for all $j\in\lra{0,\ldots,n_0-1}$, $\gamma-\alpha^{(j)}-e_{i(j)}=\gamma-\alpha^{(j+1)}=\beta^{(n_0-(j+1))}$.

Therefore, by the previous Lemma, if $\Lambda_i(\alpha^{(j)},\mc{V})\neq\vide$ then $\Lambda_i\lrp{\beta^{(n_0-(j+1))},\val(I)}=\vide$. We then obtain the result by comparing \eqref{eq:prop:cle} and \eqref{defellprime}. \hfill $\Box$

\subsubsection{End of the proof of Theorem~\ref{symmetry-values}}
\label{fin}

We are now able to finish the proof of Theorem \ref{symmetry-values}.

\s

%\textbf{End of the proof of Theorem~\ref{symmetry-values}}. 
We assume now that $D$ is Gorenstein.
 The inclusion $\val(I^\vee)\subseteq \mc{V}$ holds by Proposition~\ref{sens1}. It remains to prove that this inclusion cannot be strict.

Since $ \sum_{i=1}^p \lrp{\nu_i-\lambda_i}=\dim_{\C} t^\lambda\co_{\wt{D}}/t^\nu\co_{\wt{D}}=\dim_\C t^\lambda\co_{\wt{D}}/I+\dim_\C I/t^\nu\co_{\wt{D}}$, we have by Proposition~\ref{duality:MCM}:
$$\dim_\C I^\vee/t^{\gamma-\lambda}\co_{\wt{D}} = \sum_{i=1}^p \lrp{\nu_i-\lambda_i}\ -\dim_\C I/t^\nu\co_{\wt{D}}.$$

Thanks to Proposition~\ref{prop:cle}, we obtain:
\begin{equation}
\label{ineg}
\ell\leqslant \dim_\C I^\vee/t^{\gamma-\lambda}\co_{\wt{D}}.
 \end{equation}
However, by Lemma~\ref{lem:lp},  if $\mc{V}\neq \val(I^\vee)$, then $\ell\geqslant 1+\dim_\C I^\vee/t^{\gamma-\lambda}\co_{\wt{D}}$, which contradicts~\eqref{ineg}. \hfill $\Box$ 

%Therefore, we have $\mc{V}=\val(I^\vee)$, that is to say:  $v\in \mathrm{val}(I^\vee)\iff \Delta(\gamma-v-1,I)=\emptyset$. 

\s

Another consequence of the equality $\mc{V}=\val(I^\vee)$ is that the number $\ell$ is equal to the dimension of $I^\vee/t^{\gamma-\lambda}\co_{\wt{D}}$. Therefore, the inequality in Proposition~\ref{prop:cle} is in fact an equality. 
Moreover, since for all $w\in\Z^p$, there exist  $\lambda',\nu'\in\Z^p$ such that $\gamma-\lambda'+\N^p \subseteq \val(I^\vee)\subseteq \gamma-\nu'+\N^p$ and $\gamma-\nu'\leqslant w\leqslant \gamma-\lambda'$, it also has the following consequence:

\begin{cor} 
\label{egalite:ell:prime}
Let $D$ be a Gorenstein curve, $I\subset Q(\co_D)$ be a fractional ideal and $w\in\Z^p$. Then:
\begin{equation}
\label{equiv:lam}
\Lambda_i(w,\val(I^\vee))\neq \emptyset \iff \Lambda_i(\gamma-w-e_i,\val(I))=\emptyset.
\end{equation}
\end{cor}

\begin{cor}
\label{valeq}
Let $I, J$ be fractional ideals and $\nu\in\N^p$ be such that $\nu+\N^p\subseteq J$. We assume $\val(J)\subseteq \val(I)$. If $\dim_\C I/t^\nu\co_{\wt{D}}=\dim_\C J/t^\nu\co_{\wt{D}}$ then $\val(J)=\val(I)$.
\end{cor}
\begin{proof}
 If $\val(I)\neq \val(J)$, then as in Proposition~\ref{intrus}, there exists $w\in\val(I)\backslash \val(J)$ and $j\in\unp$ such that $\Lambda_j(w,\val(J))=\emptyset$. The same argument as in the second step of the proof of Theorem~\ref{symmetry-values} shows that $\dim_\C I/t^\nu \co_{\wt{D}}>\dim_\C J/t^\nu \co_{\wt{D}}$. Hence the result. 
\end{proof}

\subsection{Poincar\'e series of a fractional ideal}

This section follows a suggestion of Antonio Campillo. Let $(D=D_1\cup\dots\cup D_p,0)$ be the germ of a reduced reducible Gorenstein curve, and $I\subset Q(\co_D)$ be a fractional ideal. 

The following definitions are inspired by \cite{campilloandco}. We recall that $I_v=\lra{g\in I; \val(g)\geqslant v}$.

We consider the set of formal Laurent series  $\mc{L}=\Z[[t_1^{-1},\ldots,t_p^{-1},t_1,\ldots,t_p]]$ as in \cite{campilloandco}. This set is not a ring, it is only a set of formal infinite sums indexed  by $\Z^p$, with a structure of  $\Z[t_1^{-1},\ldots,t_p^{-1},t_1,\ldots,t_p]$-module. % An element in $\mc{L}$ is an infinite series of the form $\sum_{i\in \Z} a_i t^i$ with for all $i\in\Z$, $a_i\in\Z$. In particular, the multiplication of two such series is not well defined because some coefficients may be infinite, so that $\mc{L}$ is not a ring. However, the set $\mc{L}$ is a $\Z[t_1^{-1},\ldots,t_p^{-1},t_1,\ldots,t_p]$-module.

We set: 
\begin{equation}
\label{Lc}
L_I(t_1,\ldots,t_p)=\sum_{v\in\Z^p} c_I(v) t^v
\end{equation}
with $c_I(v)=\displaystyle{\dim_\C I_v/I_{v+\undun}}$ and 

\begin{equation}
\label{P}
P_I(t)=L_I(t)\prod_{i=1}^p (t_i-1).
\end{equation} 

\begin{remar}
In \cite{campilloandco}, the authors study the case $I=\co_D$ with $D$ a plane curve. They prove that $P_{\co_D}(t)$ is in fact a polynomial, and for plane curves with at least two components, $\displaystyle{\f{P_{\co_D}(t)}{t_1\cdots t_p-1}}$ is the Alexander polynomial of the curve (see \cite[ Theorem 1]{campilloandco}).
\end{remar}

Our purpose here is to deduce from Theorem~\ref{symmetry-values} a relation between $P_{I}(t)$ and $P_{I^\vee}(t)$.

The following lemma is a direct consequence of the definition  of $P_I$:
\begin{lem}
We define for $v\in\Z^p$, 
\begin{equation}
\label{alphaI}
a_I(v)=\sum_{J\se\{1,\ldots,p\}} (-1)^{\mathrm{Card}(J^c)} c_I(v-e_J)
\end{equation}
where we denote for $J=\{j_1,\cdots,j_k\}$, $e_J=e_{j_1}+\cdots+e_{j_k}$ and $J^c$ the complement of $J$ in $\lra{1,\ldots,p}$. Then  $$P_I(t)=\sum_{v\in\Z^p} a_I(v)t^v.$$
\end{lem}

We use the previous lemma to prove the following property:
\begin{lem}
The formal Laurent series $P_I(t)$ is a polynomial. 
\end{lem}
\begin{proof}
Let $\lambda,\nu\in\Z^p$ be such that $\nu+\N^p\subseteq \val(I)\subseteq \lambda+\N^p$. The only possibly non-zero $a_I(v)$ are those such that $\lambda\leqslant v\leqslant \nu$. Indeed, let us assume for example that $v_p<\lambda_p$ or $v_p>\nu_p$. We can then prove thanks to Corollary~\ref{suite} that for all $J\subset \{1,\ldots,p\}$ such that $p\notin J$, $c_I(v-e_{J\cup\{p\}})=c_I(v-e_J)$. By  definition~\eqref{alphaI}, this gives us the result. 
\end{proof}

 The symmetry of Theorem~\ref{symmetry-values} has the following consequence:

\begin{prop}
\label{sym:poly}
With the same notations, 

\begin{equation}
\label{sym:poly:mys}
P_{I^\vee}(t)=(-1)^{p+1}\ t^\gamma\ P_{I}\lrp{\f{1}{t_1},\ldots,\f{1}{t_p}}.
\end{equation}
\end{prop}
\begin{proof}
The property \eqref{sym:poly:mys} is in fact equivalent to the following property:
\begin{equation}
\label{prop:10:2}
\forall v\in\Z^p,\ a_{I^\vee}(v)=(-1)^{p+1} a_I(\gamma-v).
\end{equation}

This property is  obvious if $v\notin \lra{\w\in\Z^p;\ \gamma-\nu\leqslant w\leqslant\gamma-\lambda}$ since both $a_{I^\vee}(v)$ and $a_I(\gamma-v)$ are zero.

By \eqref{alphaI}, it is sufficient to prove that for all $v\in\Z^p$, $c_{I^\vee}(v)=p-c_I(\gamma-v-1)$. We have:
$$c_{I^\vee}(v)=\mathrm{Card}\{ i\in\{1,\ldots,p\}\ ;\ \Lambda_i(v+e_1+\cdots+e_{i-1},\val(I^\vee))\neq\vide\},$$
$$c_{I}(\gamma-v-1)=\mathrm{Card}\{ i\in\{1,\ldots,p\}\ ;\ \Lambda_i(\gamma-v-e_1-\cdots-e_i,\val(I))\neq\vide\}.$$

The result follows from the equivalence~\eqref{equiv:lam}.
\end{proof}

\section{On the structure of the set of values of \loga{} residues}

\label{part2}

In this part we give several properties of the module of logarithmic residues  along plane curves or complete intersection curves. We first recall some definitions from \cite{saitolog}. We then focus on the module of logarithmic residues along plane curves. We study the set of its zero divisors and we also give its conductor thanks to Theorem \ref{symmetry-values}. We then recall definitions from \cite{alektsikh} and \cite{alekres} for complete intersections. We prove that the values of multi-residues are in relation with the values of K\"ahler differentials. %K\"ahler differentials are used for the analytic classification of plane branches proposed in \cite{hefez}, and \cite{Hefez2} for two branches. 

\subsection{Preliminaries on \loga{} residues}

We recall here some definitions and results about \loga{} vector fields, \loga{} differential forms and their residues, which can be found in~\cite{saitolog}. In this subsection, we consider hypersurfaces, and we will focus on the case of plane curves in subsection \ref{zero divisors}.

\medskip

Let us  consider a reduced hypersurface germ $D\subset (\C^n,0)$ defined by $f\in\C\lra{x_1,\ldots,x_n}$. We denote by $\Theta_n$ the module of germs of holomorphic vector fields on $(\C^n,0)$ and $\Cx=\C\lra{x_1,\ldots,x_n}$. We set $\Omega^q_{\C^n}$, or $\Omega^q$ for short, the module of holomorphic differential $q$-forms.

\begin{de} 
\label{def}A germ of vector field $\delta\in \Theta_n$ is called logarithmic along $D$ if  $\delta(f)=\alpha f$ with $\alpha\in\Cx$. We denote by $\derl$ the $\Cx$-module of \loga{} vector fields along $D$.

 A germ of meromorphic  $q$-form $\w\in\f{1}{f}\Omega^q$ with simple poles along $D$ is called logarithmic if $f\dd\w$ is holomorphic. We denote by $\Omega^q(\log D)$ the $\Cx$-module of \loga{} $q$-forms on $D$.
\end{de}

\begin{lem}[\protect{\cite[Lemma 1.6]{saitolog}}]
The two modules $\Omega^1(\log D)$ and $\derl$ are reflexive and each is the dual $\Cx$-module of the other. 
\end{lem}
\begin{de}
If $\derl$ (or equivalently $\Omega^1(\log D)$) is a free $\Cx$-module, we call $D$ a germ of free divisor.
\end{de}

In particular, plane curves are free divisors (see \cite[1.7]{saitolog}).

\s

\begin{prop}[Saito criterion, \protect{\cite[1.8]{saitolog}}]
\label{saito:criterion}
The germ $D$ is  free \ssi{} there exists $(\delta_1,\ldots,\delta_n)$ in $\der(-\log D)$ such that $\delta_j=\sum a_{ij} \dr_{x_i}$ with $\det\lrp{(a_{ij})_{1\leqslant i,j\leqslant n}}=uf$, where $u$ is invertible in $\Cx$. 
\end{prop} 

In order to define the notion of \loga{} residues, we need the following characterization of \loga{} differential forms:

\begin{prop}[\protect{\cite[1.1]{saitolog}}]
\label{prop:loga:form}
A meromorphic $q$-form $\w$ with simple poles along $D$ is \loga{} \ssi{} there exist $g\in\Cx$, which does not induce a zero divisor in  $\co_D=\Cx/(f)$, a holomorphic $(q-1)$-form  $\xi$ and a holomorphic $q$-form $\eta$ such that: 
\begin{equation}
\label{forme:loga}
g\w= \f{\dd f}{f}\wedge \xi+\eta.
\end{equation}
\end{prop}

\begin{de}
The residue $\res^{q}(\w)$ of $\w\in \Omega^q(\log D)$ is defined by $$\res^q(\w):=\restriction{\f{\xi}{g}}{D}\in Q(\co_D)\otimes_{\co_D} \Omega^{q-1}_D,$$ where $\xi$ and $g$ are given by~\eqref{forme:loga}, and $\displaystyle{\Omega^{q-1}_D=\restriction{\frac{\Omega_{\C^n}^{q-1}}{\dd f\wedge \Omega_{\C^n}^{q-2}+f\Omega^{q-1}_{\C^n}}}{D}}$ is the module of K\"ahler differentials on $D$.
\end{de}

If $q=1$, we write $\res(\w)$ instead of $\res^1(\w)$. 
\begin{de}
We define $$\mc{R}_D:=\lra{\res(\w); \w\in\Omega^1(\log D)}\subseteq Q(\co_D).$$ We call $\mc{R}_D$ the module of \loga{} residues of $D$. 
\end{de}

In particular, $\mc{R}_D$ is a finite type $\co_D$-module. Moreover, it satisfies the following property:
\begin{lem}[\protect{\cite[Lemma 2.8]{saitolog}}]
 We have the inclusion $\co_{\wt{D}}\subseteq \mc{R}_D$
\end{lem}

\begin{nota}
We denote by $\mc{J}_D\subseteq \co_D$ the Jacobian ideal of $D$, that is to say the ideal of $\co_D$ generated by the partial derivatives of $f$. 
\end{nota}
The following result gives the relation between the module of logarithmic residues and the Jacobian ideal: 

\begin{prop}[\protect{\cite[Proposition 3.4]{gsres}}]
\label{th:dual:jd}
Let $D$ be the germ of a reduced divisor. Then $\mc{J}_D^\vee=\mc{R}_D$. If moreover $D$ is free, $\mc{R}_D^\vee=\mc{J}_D$.
\end{prop}

\subsection{Logarithmic residues along plane curves}
\label{zero divisors}

We give here several properties of the set of values of \loga{} residues $\mc{R}_D$ of a plane curve $D$. We first determine the zero divisors included in $\mc{R}_D$, see Proposition \ref{zerodiv:res}. We also relate the conductor of $\mc{R}_D$ to the multiplicities of the branches of $D$ (see Proposition~\ref{cond:res}). We end this subsection with the relation between logarithmic residues and the torsion of K\"ahler differentials.

%\subsubsection{Generators}

Since plane curves are free divisors, the module $\Omega^1(\log D)$ is a free $\Cxy$-module of rank $2$. Let us assume that $\w_i=\frac{\alpha_i\dd_x+\beta_i\dd_y}{f},\ i=1,2$ is a basis of $\Omega^1(\log D)$.
% Then by duality, a basis of $\Omega^1(\log D)$ is:
%\begin{align*}
%\w_1&=\f{\beta_2 dx-\alpha_2dy}{f}\\
%\w_2&=\f{-\beta_1 dx+\alpha_1 dy}{f}\\
%\end{align*}

If for $c_1,c_2\in\C$, $g=c_1\cdot f_x'+c_2\cdot f_y'$ induces a non zero divisor in $\co_D$, then the module of residues is generated by
$\displaystyle{\res(\w_1)=\frac{c_1\cdot \alpha_1+c_2\cdot \beta_1}{g}}$ and $\displaystyle{
\res(\w_2)=\frac{c_1 \cdot \alpha_2+c_2 \cdot \beta_2}{g}}$. Thus, the module of logarithmic residues can be generated by two elements.

\subsubsection{Zero divisors}

Let $D=D_1\cup\dots\cup D_p$ be the germ of a reduced plane curve defined by $f=f_1\cdots f_p$ where for all $i\in\unp$, $f_i$ is irreducible.

We first want to prove that the negative values of $\mc{R}_D$ determine all the values of $\mc{R}_D$. It comes from a general property of fractional ideals. %Since we will need a similar result for the Jacobian ideal, we state it for an arbitrary fractional ideal. 

We recall that for $g\in Q(\co_D)$, $\val_i(g)=\infty$ means that the restriction of $g$ on $D_i$ is zero.

The following proposition shows that the values of the zero divisors are determined by the faces of the negative quadrant with origin $\nu$.
\begin{prop}
\label{zero:div}
Let $I\subset Q(\co_D)$ be a fractional ideal and let $\nu\in\Z^p$ be any element satisfying $t^\nu\co_{\wt{D}}\subseteq I$. Let $\alpha\in\lrp{\Z\cup\lra{\infty}}^p$. Then $\alpha\in\ov{\val(I)}$ if and only if either $\alpha\in\val(I)$, or the element $w$ defined by $w_i=\alpha_i$ if $\alpha_i\in\N$ and $w_i=\nu_i$ if $\alpha_i=\infty$ satisfies $w\in\val(I)$. 
\end{prop} 

\begin{proof}
Let $\alpha\in\ov{\val(I)}$ be such that $\alpha\notin \val(I)$. Let $\alpha'\in \nu+\N^p\subseteq \val(I)$ be such that if $\alpha_i\in\N$, then $\alpha_i'\geqslant \alpha_i$, and if $\alpha_i=\infty$, $\alpha_i'=\nu_i$. Then, by Proposition~2.9, $\inf(\alpha,\alpha')=w\in\val(I)$ with $w_i=\alpha_i$ if $\alpha_i\in\N$ and $w_i=\nu_i$ if $\alpha_i=\infty$. 

Conversely, let $w\in \val(I)$. Let us assume that there exists $j\in\lra{1,\ldots,p}$ such that $w_j=\nu_j$.  Let $J$ be a set of indices such that for all $j\in J$, $w_j=\nu_j$. Let us show that the element $\alpha$ defined by $\alpha_i=w_i$ if $i\notin J$ and $\alpha_i=\infty$ if $i\in J$ satisfies $\alpha\in\ov{\val(I)}$.  Let $h\in I$ be such that $\val(h)=w$. Since $t^\nu\co_{\wt{D}}\subseteq I$, there exists $g\in I$ such that for all $j\in J$, $g|_{D_j}=h|_{D_j}$ and for all $j\notin J$, $\val_j(g)>w_j$. Then $h-g$ is a zero divisor of $I$ whose value $\alpha$ satisfies for all $j\in J$, $\alpha_j=\infty$ and for all $j\notin J$, $\alpha_j=w_j$. 
\end{proof}
 
%\begin{remar}
We will use the notation $\mc{M}_I=\lra{w\in\val(I)\ ;\ \mathrm{Card}\lra{j\in\lra{1,\ldots,p}\ ;\ w_j=\nu_j}\geqslant 1}$. This set determines the value of the zero divisors of $Q(\co_D)$. 
%\end{remar} 
%\begin{prop}
%\label{zero:div}
%Let $I\subset Q(\co_D)$ be a fractional ideal and let $\nu\in\Z^p$ be such that $t^\nu\co_{\wt{D}}\subseteq I$.  The values of the zero divisors of $I$ are determined by $$\mc{M}_I=\bigcup_{q=1}^{p} \bigcup_{\sigma\in \mathfrak{S}_q} \lra{w\in\val(I); \forall i\in\lra{1,\ldots,q}, w_{\sigma(i)}=\nu_{\sigma(i)}}$$
%where $\mathfrak{S}_q$ is the group of permutations of $q$ elements. More precisely, the values of the zero divisors are exactly the values obtained by replacing for $v\in\mc{M}_I$ the coordinates which satisfies $v_j=\nu_j$ by $\infty$. 
%\end{prop}
%\begin{proof}
%Let $g\in I$ be a zero divisor. Then by Proposition \ref{propinf}, $\inf(\val(g),\nu)\in\val(I)$. It is then obvious that $\inf(\val(g),\nu)\in\mc{M}_I$.
%
%Conversely, let $w\in\mc{M}_I$.  Let $h\in I$ be such that $\val(h)=w$. Let $\emptyset \neq J\subseteq \unp$ be a set of indices such that for all $j\in J$, $w_j=\nu_j$. Since $t^\nu\co_{\wt{D}}\subseteq I$, there exists $g\in I$ satisfying  for $j\in J$, $g|_{D_{j}}=h|_{D_{j}}$, and  $\val_j(g)>w_j$ if $j\notin J$. Then $h-g$ is a zero divisor whose value satisfies for all $j\in J$, $\val_{j}(h-g)=\infty$ and for all $j\notin J$, $\val_j(h-g)=w_j$. 
%\end{proof}

\begin{cor}
\label{cadran}
Let $I\subset Q(\co_D)$ be a fractional ideal and $\nu\in\mathbb{Z}^p$ be such that $\nu+\N^p\subseteq\val(I)$. Let $w\in\Z^p$. Then $$w\in\val(I) \iff \inf(w,\nu)\in\val(I).$$ In particular, it means that the set $$\val(I)\cap \lra{w\in\Z^p; w\leqslant \nu}$$ determines the set $\val(I)$. 
\end{cor}
\begin{proof}
The implication $\Rightarrow$ comes from Proposition~\ref{propinf}. 
For the implication $\Leftarrow$, let $w\in\Z^p$ be such that $\inf(w,\nu)\in\val(I)$.  If $w\leqslant \nu$, then  $w=\inf(w,\nu)\in\val(I)$. If there exists $j$ such that $w_j>\nu_j$, then $\inf(w,\nu)\in\mc{M}_I$ where $\mc{M}_I$ is defined after Proposition~\ref{zero:div}. By Proposition~\ref{zero:div}, there exist a zero divisor $g\in I$ such that $\val_j(g)=\infty$ if $w_j\geqslant \nu_j$ and $\val_j(g)=w_j$ if $w_j<\nu_j$.   Let $v=\max(w,\nu)$. Since $\nu+\N^p\subseteq \val(I)$, we have $v\in\val(I)$ and $w=\inf(v,\val(g))\in \val(I)$.
\end{proof}
%Therefore, since $t^\nu\co_{\wt{D}}\subseteq I$, there exists a zero divisor $g\in I$, $q\in \unp$ and $\sigma\in\mathfrak{S}_q$ which satisfy for all $i\in\lra{1,\ldots,q}$, $\val_{\sigma(i)}(g)=\infty$ if $v_{\sigma(i)}=\nu_{\sigma(i)}$ and for all $j\notin \sigma\lrp{\lra{1,\ldots,q}}$, $\val_j(g)=v_j<\nu_i$. 
\begin{remar}
By Proposition \ref{zero:div} and Corollary~\ref{cadran}, the set $\val(I)\cap \lra{w\in\Z^p; w\leqslant \nu}$ also determines $\ov{\val(I)}$.
\end{remar}

The inclusion $\co_{\wt{D}}\subseteq \mc{R}_D$ gives the following corollary:

\begin{cor}
\label{neg:res}
The set of values of $\mc{R}_D$ is determined by the set $$\lra{v\in\val(\mc{R}_D); v\leqslant 0}.$$ More precisely, we have:

$$\val(\mc{R}_D)=\lra{v\in\val(\mc{R}_D); v\leqslant 0}\cup \lra{v\in\Z^p; \inf(v,0)\in\val(\mc{R}_D)}.$$
\end{cor}

Let us determine the values of $\mc{R}_D$ which come from the branches or union of branches. 

%$q\in\lra{1,\ldots,p}$ and $D'=D_{i_1}\cup\dots\cup D_{i_q}$ be the union of $q$ branches of $D$.Moreover, $v\in\lra{\alpha\in\val(\mc{R}_D); \forall j\notin\lra{1,\ldots,q}, \alpha_{i_j}=0}$ if and only if there exists a logarithmic $1$-form $\w\in \Omega^1\lrp{\log D'}$ such that  for all $j\notin \lra{i_1,\ldots, i_q}$, $\val_j(\res(\w))=v_j$, and for all $j\in\lra{1,\ldots,q}$, $\val_{i_j}(\res(\w))=\infty$. 
\begin{prop}
\label{zerodiv:res}
Let $\emptyset\neq J\subseteq \lra{1,\ldots,p}$ and $D'=\bigcup_{j\in J} D_j$. Then $\Omega^1(\log D')\subseteq \Omega^1(\log D)$.

Renumbering the branches, we may assume that $J=\lra{1,\ldots, q}$ with $q\leqslant p$. Then: 
$$\mc{R}_D\cap \lrp{Q(\co_{D_1})\oplus \dots \oplus Q(\co_{D_q})\oplus\lra{0}^{p-q}}=\mc{R}_{D'}.$$ 
\end{prop}

%\begin{prop}
%\label{zerodiv:res}
%Let $\emptyset\neq J\subseteq \unp$, and $D'=\bigcup_{j\in J} D_j$.  Then $$\Omega^1(\log D')\subseteq \Omega^1(\log D)$$
%
%Moreover, $v\in\lra{\alpha\in\val(\mc{R}_D); \forall j\notin J, \alpha_j=0}$ if and only if there exists a logarithmic $1$-form $\w\in \Omega^1\lrp{\log D'}$ such that  for all $j\in J$, $\val_j(\res(\w))=v_j$. 
%\end{prop}
\begin{proof}
For the first part of the statement,  we set $F_1$ the equation of $D'$.  Let $\w$ be a \loga{} $1$-form along $D'$. Then, $F_1\w$ and $F_1 \dd\w$ are holomorphic, so that $f\w$ and $f \dd\w$ are holomorphic. 

For the second part of the statement, let us notice the following property. Let $\w\in\Omega^1(\log D)$. Then $\w\in\Omega^1(\log D')$ if and only if for all $j\notin J$, $\val_j(\res(\w))=\infty$. The second part of the statement then follows from  this remark and Proposition~\ref{zero:div}.
\end{proof}

In particular, the logarithmic residues of the irreducible components satisfy the following property: 

\begin{cor}
\label{somme:directe}
We have the following inclusion:
$$\mc{R}_{D_1}\oplus\dots\oplus \mc{R}_{D_p}\incl \mc{R}_D.$$
Therefore, $\val_1(\mc{R}_{D_1})\times\dots\times\val_p(\mc{R}_{D_p})\subseteq \val(\mc{R}_D)$.
\end{cor}

\begin{remar}
If $D=D_1\cup D_2$ is a plane curve satisfying $\mc{R}_D=\mc{R}_{D_1}\oplus\mc{R}_{D_2}$, then by \cite{snc}, it is a splayed divisor, and in fact it is even a normal crossing plane curve. We refer to \cite[Definition 2.3]{fabersplayed} for the definition of a splayed divisor. 
\end{remar}

%\begin{remar}
%Let us recall that a divisor $D\subset (\C^n,0)$ is called \emph{splayed} if there exist local coordinates $(y_1,\ldots,y_n)$ in a neighbourhood of $0\in\C^n$, $1\leqslant p< n$, $g_1\in\C\lra{y_1,\ldots,y_p}$, $g_2\in\C\lra{y_{p+1},\ldots,y_n}$ such that the defining equation of $D$ is $f(y_1,\ldots,y_n)=g_1(y_1,\ldots,y_p)\cdot g_2(y_{p+1},\ldots,y_n)$ (see for example \cite{fabersplayed} for more details on splayed divisors). If $D=D_1\cup D_2$ is a plane curve satisfying $\mc{R}_D=\mc{R}_{D_1}\oplus\mc{R}_{D_2}$, then by \cite{snc}, it is a splayed divisor, and in fact it is even a normal crossing plane curve. 
%\end{remar}

We now study the set of values of the dual of $\mc{R}_D$, namely,  the Jacobian ideal $\mc{J}_D$. We show that the modules of \loga{} vector fields $\der(-\log D_i)$ for $i\in\unp$ give information on the structure of the set of values of the Jacobian ideal.

%$j\notin\lra{i_1,\ldots,i_q}$, $\val_j(\delta(f_j))=v_j-\sum_{i\neq j} \val_j(f_i)$. 
\begin{prop}
\label{zero:div:jac}
Let $\emptyset\neq J\subseteq \lra{1,\ldots,p}$ and $D'=\bigcup_{j\in J} D_j$. 
Renumbering the branches, we may assume that $J=\lra{1,\ldots, q}$ with $q\leqslant p$. Then:
$$\mc{J}_D\cap \lrp{\lra{0}^{q}\oplus \C\lra{t_{q+1}}\oplus\dots\oplus \C\lra{t_p}}=\lra{\delta(h)\ ;\ \delta\in\der(-\log D')}.$$  In particular, the set of zero divisors of $\mc{J}_D$ is determined by the family of modules $$\lra{\mathrm{Der}\big(-\log (\cup_{j\in J} D_j)\big)}_{J\subset \unp, J\neq \emptyset}.$$ 
\end{prop}
%\begin{prop}
%\label{zero:div:jac}
%Let $v\in\Z^p$. We keep  the notations of Proposition~\ref{zero:div}. Then  $v\in\mc{M}_{\mc{J}_D}$ if and only if there exists $D'=\bigcup_{j\in J} D_j$ with $\emptyset \neq J\subseteq \unp$ and $\delta\in\der(-\log D')=\bigcap_{j\in J} \mathrm{Der}(-\log D_j)$ such that for $j\notin J$, $\val_j(\delta(f_j))=v_j-\sum_{i \neq j} \val_j(f_i)$.  In particular, the set of zero divisors of $\mc{J}_D$ is determined by the family of modules $\bigcap_{j\in J} \mathrm{Der}(-\log D_j)$ for $J\subset \unp$.
%\end{prop}

%, up to a permutation of indices: $$\val(\delta(f))=(\infty,\ldots,\infty,v_k,\ldots,v_p)$$ which is equivalent to:  $$\forall i\in\{1,\ldots,k-1\}, \delta(f_i)\in (f_i) \text{ and } \forall i\in \{k,\ldots,p\}, \val_i(\delta(f))=v_i$$Hence the result.
\begin{proof}
We first notice that for all $g\in\mc{J}_D$, there exists $\delta\in\Theta_2$ such that $\delta(f)=g$ in $\co_D$, where $\Theta_2$ is the module of holomorphic vector fields on $(\C^2,0)$. Moreover, $\delta(f)$ induces in $\co_{\wt{D}}=\prod_{i=1}^p \co_{\wt{D_i}}$ the element: $$\delta(f)=\lrp{f_2\cdots f_p\delta(f_1),\ldots,f_1\cdots f_{p-1}\delta(f_p)}.$$ By Proposition~\ref{zero:div}, $v\in\mc{M}_{\mc{J}_D}$ if and only if there are $\emptyset\neq J\subseteq \unp$ and $\delta\in \Theta_2$ such that for all $j\in J$, $\val_j(\delta(f))=\infty$ 	and for all $j\notin J$, $\val_j(\delta(f))=v_j$. This condition is equivalent to the condition:  for all $j\in J$, $\delta(f_j)\in (f_j)$ and  for all $j\notin J$, $\val_j(\delta(f_j))=v_j-\sum_{i \neq j} \val_j(f_i)$.
\end{proof}

\subsubsection{Conductor of the module of residues}

We compute here the conductor of $\mc{R}_D$, that is to say, the minimal $\nu\in\Z^p$ such that $\nu+\N^p\subseteq \val(\mc{R}_D)$.  

We need the following results:
\begin{prop}[\protect{\cite[Theorem 2.7]{delgado2}}]
\label{cond:delgado}
Let $f=f_1\cdots f_p$ be a reduced equation of a plane curve germ. We assume that for all $i\in\unp$, $f_i$ is irreducible. We denote by $c_i$ the conductor of the branch defined by $f_i$. The conductor of $D$ is given by $$\gamma=\lrp{c_1+\sum_{i=2}^p {\mathrm{val}_1(f_i)},\ldots,c_p+\sum_{i=1}^{p-1} {\mathrm{val}_p(f_i)}}.$$
\end{prop} 

We then have:

\begin{lem}
\label{val:jc:xy}
Let $f\in\Cxy$ be a reduced equation of a plane curve germ. Then:
\begin{align*}
\val(f_x')&=\gamma+\val(y)-\undun\\
\val(f_y')&=\gamma+\val(x)-\undun
\end{align*}
\end{lem}
\begin{proof}
If $f$ is irreducible, it is exactly the statement of Teissier's lemma (see~\cite[Lemma 2.3]{cassouploski}). If $f$ is reducible, we prove the result for $f_x'$. We have the following  equality: $$\val_j\lrp{f_x'}=\sum_{i\neq j}\val_j(f_i) +\val_j\lrp{\f{\dr f_j}{\dr x}}.$$ By  Teissier's Lemma, $\val_j\lrp{f_x'}=c_j+\val_j(y)-1$.  Theorem~\ref{cond:delgado} then gives the result.
\end{proof}

\begin{prop}
\label{cond:res}
The conductor of $\mc{R}_D$ is $-(m^{(1)},\ldots,m^{(p)})+\undun$, where $m^{(j)}=\inf(\val_j(x),\val_j(y))$ is the multiplicity of the component $D_j$ of $D$. 
\end{prop}
\begin{proof}
By Lemma~\ref{val:jc:xy}, $\inf(\val(\mc{J}_D))=\gamma+\inf(\val(x),\val(y))-\undun$. Therefore, $$\sup\lra{\alpha\in\Z^p\ ; \forall \beta\leqslant \alpha, \Delta(\beta,\mc{J}_D)=\emptyset}=\gamma+\inf(\val(x),\val(y))-\underline{2}.$$ By Theorem~\ref{symmetry-values}, the conductor of $\mc{R}_D$ is $\nu=-\inf(\val(x),\val(y))+\undun$.  
\end{proof}

\subsubsection{Relation with the torsion of K\"ahler differentials}

We mention here the relation between \loga{}  forms and the torsion of K\"ahler differentials for a  plane curve $D$. It leads to a determination of the dimension  of the torsion  of $\Omega^1_D$ as a $\C$-vector space when $D$ is a plane curve which is slightly different from the proofs of O.Zariski (see \cite{zariski-torsion}) and R.Michler (see \cite{michler-torsion}).

We first assume that $D$ is the germ of a reduced hypersurface in $(\C^n,0)$. The following property was proved by A.G.Aleksandrov:

\begin{prop}[\protect{\cite[3.1]{alektorsion}}]
\label{prop:alek:torsion}
For all $1\leqslant q\leqslant n$, the following map:

\begin{align*}
\frac{\Omega^{q}(\log D)}{\frac{\dd f}{f}\wedge \Omega^{q-1}_{\C^n}+\Omega^{q}_{\C^n}} &\to \mathrm{Tors}(\Omega^{q}_D)\\
 [\w]&\mapsto [f\w]
 \end{align*}
 is an isomorphism of $\co_D$-modules.
\end{prop}

\begin{proof}
It is a consequence of the characterization~\eqref{forme:loga} of logarithmic forms.
\end{proof}

\begin{cor}
The map $\res^{q}$ induces an isomorphism of $\co_D$-modules: 
$$\frac{\res^{q}(\Omega^{q}(\log D))}{\Omega^{q-1}_D} \simeq \mathrm{Tors}(\Omega^{q}_D).$$
\end{cor} 
\begin{proof}
We have the following exact sequences:

$$0\to\Omega^q_{\C^n}\to \Omega^q(\log D)\to \res^q\lrp{\Omega^q(\log D)}\to 0,$$
$$0\to \Omega^q_{\C^n}\to \frac{\dd f}{f}\wedge \Omega^{q-1}_{\C^n}+\Omega^q_{\C^n}\xrightarrow{\res^q} \Omega^{q-1}_D\to 0.$$
Therefore, $\res^{q}(\Omega^{q}(\log D))\simeq \frac{\Omega^q(\log D)}{\Omega^q_{\C^n}}$ and $\Omega^{q-1}_D\simeq \frac{ \frac{\dd f}{f}\wedge\Omega^{q-1}{\C^n}+\Omega^q_{\C^n}}{\Omega^q_{\C^n}}$, so that by a classical isomorphism theorem, we have $\frac{\res^{q}(\Omega^{q}(\log D))}{\Omega^{q-1}_D} \simeq\frac{\Omega^{q}(\log D)}{\frac{\dd f}{f}\wedge \Omega^{q-1}_{\C^n}+\Omega^{q}_{\C^n}}$, and we conclude with Proposition~\ref{prop:alek:torsion}.
\end{proof}

\begin{cor}
Let $D\subseteq(\C^2,0)$ be a plane curve germ. We denote by $\tau=\dim_C \co_D/\mc{J}_D$ the Tjurina number. Then $\displaystyle{\f{\mc{R}_D}{\co_D}\simeq \mathrm{Tors}(\Omega^{1}_D)}$ and $\dim_{\C} \mathrm{Tors}(\Omega^{1}_D)=\tau$. 
\end{cor}
\begin{proof}
We use Propositions~\ref{duality:MCM} and \ref{th:dual:jd} to prove that $\dim_\C \mc{R}_D/\mc{O}_D=\dim_{\C} \mc{O}_D/\mc{J}_D=\tau$.
\end{proof}

\subsection{Complete intersection curves and the relation with K\"ahler differentials}

\label{complete-int}

This section is devoted to the study of complete intersection curves, which are a particular case of Gorenstein curves. 

We begin with the definition of multi-logarithmic forms  along a reduced complete intersection given in \cite{alekres}. We then focus on the case of complete intersection curves for which we give the relation between the values of multi-residues and the values of K\"ahler differentials.

\s

Let $C\subset (\C^m,0)$ be the germ of a  reduced complete intersection defined by a regular sequence $(h_1,\ldots,h_k)$.  For $j\in\lra{1,\ldots,k}$, we set $\wh{h_j}=h_1\cdots h_{j-1}\cdot h_{j+1}\cdots h_k$.

\begin{de}[\cite{alekres}]
Let $\w\in \frac{1}{h_1\dots h_k}\Omega^q$ with $q\in \mathbb{N}$. Then $\w$ is called a \emph{multi-logarithmic differential $q$-form} along the complete intersection $C$ if 

$$\forall i\in\unk,\  dh_i\wedge \w\in\sum_{j=1}^{k} \f{1}{\wh{h_j}}\Omega^{q+1}.$$ 

We denote by $\Omega^q(\log C)$ the $\Cx$-module of multi-logarithmic $q$-forms along $C$.
\end{de}
To simplify the notations, we set $\wt{\Omega}^q:=\displaystyle{\sum_{j=1}^{k} \f{1}{\wh{h_j}}\Omega^{q}}$.

If $k=1$, the definition of multi-\loga{} forms coincides with the definition of \loga{} forms \ref{def}. 

Then we have the following characterization which should be compared with Proposition \ref{prop:loga:form}: 

\begin{theo}[\protect{\cite[\S 3, Theorem 1]{alekres}}]
\label{theo:alek}
Let $\w\in\f{1}{h_1\cdots h_k}\Omega^q$, with $q\geqslant k$. Then $\w\in\Omega^q(\log C)$ if and only if  there exist a holomorphic function $g\in\Cx$ which does not induce a zero divisor in $\co_C$, a holomorphic differential form $\xi\in\Omega^{q-k}$ and a meromorphic $q$-form $\eta\in\wt{\Omega}^q$ such that:

\begin{equation}
\label{def:k:log}
g\w=\xi\wedge\f{\dd h_1\wedge\dots\wedge \dd h_k}{h_1\cdots h_k}+\eta.
\end{equation}
\end{theo}

We set for $q\geqslant 0$ :
$$\Omega^{q}_C=\restriction{\f{\Omega^{q}}{(h_1,\ldots,h_k)\Omega^{q}+ \dd h_1\wedge \Omega^{q-1}+\cdots+ \dd h_k\wedge \Omega^{q-1}}}{C}.$$

\begin{de}
Let $\w\in\Omega^q(\log C)$, $q\geqslant k$. Let us assume that $g,\xi,\eta$ satisfy the properties of Theorem~\ref{theo:alek}. Then the  \emph{multi-residue}  of $\w$ is: 

$$\res^q_C(\w):=\restriction{\f{\xi}{g}}{C}\in Q(\co_C) \otimes_{\co_{C}} \Omega^{q-k}_C=Q(\co_{\wt{C}})\otimes_{\co_{\wt{C}}}\Omega^{q-k}_{\wt{C}}.$$

We define $\mc{R}_C^{q-k}:=\res^q_C(\Omega^{q}(\log C))$. 
In particular, if $q=k$, $\res^k_C(\w)\in Q(\co_C)$, and we denote $\mc{R}_C:=\res^k_C\left(\Omega^k(\log C)\right)$.
\end{de}

It is proved in \cite{alekres} that for $\w\in\Omega^q(\log C)$ the multi-residue $\res^q_C(\w)$ is well-defined with respect to the choices of $\xi, g$ and $\eta$ in \eqref{def:k:log}.

\s

\begin{prop}[\protect{\cite[Lemma 5.4]{snc}}]
\label{dual:snc}
Let $\mc{J}_C\subseteq \co_C$ be the Jacobian ideal, that  is to say the ideal of $\co_C$ generated by the $k\times k$ minors of the Jacobian matrix. Then:
$$\mc{J}_C^\vee=\mc{R}_C.$$
\end{prop}

\begin{remar}
 In \cite{polfreeci}, we give a more direct proof of this duality, which is not based on the isomorphism between the module of  multi-residues and the module of regular meromorphic forms given in \cite[Theorem 3.1]{alektsikh}.
\end{remar}

From now on, we assume that $C=C_1\cup \dots\cup C_p$ is a reduced complete intersection curve defined by a regular sequence $(h_1,\ldots, h_{m-1})$. We denote by $\varphi_i(t_i)=(x_{i,1}(t_i),\ldots,x_{i,m}(t_i))$ a parametrization of the branch $C_i$, which is induced by a normalization of $C$.

\begin{de}
\label{kahler:val}
Let $\w=\sum_{j=1}^m a_j \dd x_j\in\Omega^1_C$. We set $x_{i,j}'$ for the derivative of $x_{i,j}$ with respect to $t_i$ and $\varphi_i^*(\w)=\lrp{\sum_{j=1}^m a_j\circ \varphi_i(t_i) x_{i,j}'(t_i)}\dd t_i$. Then: $$\val_i(\w)=\val_i(\varphi_i^*(\w))=1+\val_i\lrp{\sum_{j=1}^m (a_j\circ\varphi_i)(t_i)\cdot x_{i,j}'(t_i)}.$$
\end{de}

We recall that $\mc{C}_C$ denotes the conductor ideal of $C$. We set $ \frac{\varphi^*(\Omega^1_C)}{\dd \underline{t}}\subset \co_{\wt{C}}=\bigoplus_{i=1}^p \C\lra{t_i}$ the fractional ideal generated by $\big( (x_{1,1}'(t_1),\ldots,x_{p,1}'(t_p)\big),\ldots, \big( (x_{m,1}'(t_1),\ldots,x_{m,p}'(t_p)\big)$. We then have $\val(\Omega^1_C)=\val\lrp{\frac{\varphi^*(\Omega^1_C)}{\dd \underline{t}}}+\undun$.

\begin{prop}
\label{gamma} 
\label{jc:kahler}
Let $C=C_1\cup \dots \cup C_p \subset \C^m$ be a reduced complete intersection curve defined by a regular sequence $(h_1,\ldots,h_{m-1})$. Then there exists $g\in\mc{C}_C$ with $\val(g)=\gamma$ such that $\mc{J}_C=g\cdot \frac{\varphi^*(\Omega^1_C)}{\dd \underline{t}}$. In particular, \begin{equation}
\label{val:jc:kah}
\val(\mc{J}_C)=\gamma+\val(\Omega_C^1)-\underline{1}.
\end{equation}
\end{prop}

\begin{proof}

Let $i\in\lra{1,\ldots,p}$.  Let $\mathrm{Jac}(h_1,\ldots,h_{m-1})$ be the Jacobian matrix of $(h_1,\ldots,h_k)$. Let $J_i$ denote the $k\times k$ minor of $\mathrm{Jac}(h_1,\ldots,h_{m-1})$ obtained by removing the column $i$. 
 Then, for all $j\in\{1,\ldots,m-1\}$ we have $h_j\circ \varphi_i(t_i)=0$, thus:

$$\Big(\mathrm{Jac}(h_1,\ldots, h_k)\circ\varphi_i(t_i)\Big)\big(x_{i,1}'(t_i),\ldots,x_{i,m}'(t_i)\big)^{t}=\big(0,\ldots,0\big)^{t}.$$

We multiply on the left by the adjoint of the matrix obtained by removing the last column of $\mathrm{Jac}(h_1,\ldots,h_k)\circ\varphi_i(t_i)$, which gives the needed relations: for all $j\in\lra{1,\ldots,m-1}$, 

$$\big(J_m\circ\varphi_i(t_i)\big)\cdot x_{i,j}'(t_i)+(-1)^{m-(j-1)} \big(J_j\circ\varphi_i(t_i)\big)\cdot x_{i,m}'(t_i) =0.$$

 We assume for example that $x_{i,m}'(t_i)\neq 0$. 
 
By setting $\displaystyle{g_i(t_i)=\frac{J_m\circ\varphi_i(t_i)}{x_{i,m}'(t_i)}}$ one obtains 
 for all $\ell\in\lra{1,\ldots,m}$,
\begin{equation}
\label{lambdai}
g_i(t_i)\cdot x_{i,\ell}'(t_i)=(-1)^{m-\ell} J_\ell\circ\varphi_i(t_i).
\end{equation}

It remains to prove that $\val_i(g)=\gamma_i$.

Let us denote by $\Pi_C$ the ramification ideal of the curve $C$, which is the $\co_{\wt{C}}$-module generated  by $$\lrp{x_{i,1}'(t_1),\ldots,x_{i,p}'(t_p)}_{1\leqslant i\leqslant m}.$$

By \cite[Corollary 1, Proposition 1]{piene}, one has: 

$$\mc{C}_C\Pi_C=\mc{J}_C\co_{\wt{C}}.$$

Thus, we have the equality $\inf(\val(\Pi_C))+\gamma=\inf(\val(\mc{J}_C))$.

The equalities~\eqref{lambdai} imply that for all $i\in\lra{1,\ldots,p}$, and $j\in\lra{1,\ldots,m}$ we have $\val_i(J_j)=\inf(\val_i(\mc{J}_C))$ if and only if $\val_i(x_{i,j}')=\inf(\val_i(\Pi_C))$.

Therefore, if $j$ is such that $\val_i(J_j)=\inf(\val_i(\mc{J}_C))$, then $\val_i(J_j)=\gamma_i+\val_i(x_{i,j}')$, which gives us $\val_i(g_i)=\gamma_i$.
\end{proof}

\begin{cor}
\label{res:kahler}
Let $C\subset \C^m$ be a reduced complete intersection curve. With the notation of Part~\ref{part1}, for all $v\in\Z^p$, we have the following equivalence:

$$v\in\val(\mc{R}_C)\iff \Delta(-v,\val(\Omega^1_C))=\emptyset$$

and $\displaystyle{\mc{R}_C=\frac{1}{g}\cdot\lrp{\frac{\varphi^*(\Omega^1_C)}{\underline{\dd t}}}^\vee}$, where $g$ is given by Proposition~\ref{gamma}.
\end{cor}
\begin{proof}
It is a consequence of the symmetry Theorem \ref{symmetry-values} together with Proposition~\ref{dual:snc}, and Proposition~\ref{gamma}
\end{proof}

\begin{remar}
\label{reg:mer:form}
The latter corollary gives also the relation between meromorphic regular forms as defined in \cite{barletdual} and K\"ahler differentials. Indeed, by \cite[Th\'eor\`eme 2.4]{alektsikh}, the module $\mc{R}_C$ of multi-residues is isomorphic to the module of regular meromorphic forms $\w_C$, which can be defined as $\w_C=\mathrm{Ext}^{m-1}_{\co_{\C^m}}\lrp{\co_C,\Omega^m_{\C^m}}$.  In particular, $\displaystyle{\w_C\simeq \frac{1}{g}\cdot\lrp{\frac{\varphi^*(\Omega^1_C)}{\underline{\dd t}}}^\vee}$.
\end{remar}

\begin{remar}
\label{jac:semi}
Another consequence of Proposition~\ref{jc:kahler} is the following inclusion:
$$\gamma+\big(\val(\co_C)\backslash\{\underline{0}\}\big)-\undun\subseteq \val(\mc{J}_C).$$

Indeed, if $h\in\m$, with $\m$ the maximal ideal of $\co_C$,  then $\val(\dd h)=\val(h)$, which gives us the inclusion $\val(\co_C)\backslash\lra{0} \subseteq \val(\Omega^1_C)$.
\end{remar}

\section{Equisingular deformations of plane curves and the stratification by \loga{} residues}

\label{part3}

The purpose of this last section is to study the behaviour of the values of \loga{} residues in an equisingular deformation of a plane curve. The base space of this deformation is denoted by $S$ in forthcoming definition~\ref{def:equi}. Each $s\in S$ is associated with a germ of plane curve $D_s$, and $\val(\co_{D_s})$ does not depend on $s$.  By the results of section~\ref{part2}, partitioning  $S$ by the values of \loga{} residues  is the same as  by the values of K\"ahler differentials. This partition is an essential ingredient of the analytic classification of plane curves described in \cite{hefez} and \cite{Hefez2} respectively for irreducible curves and for reducible curves with two branches. %They prove that in a stratum of the stratification by the values of K\"ahler differentials, the analytic classification can be expressed in terms of a normal parametrization of the curve.  The analytic equivalence is then represented by the group action of the roots of the unity of a certain order determined by $\val(\Omega^1_D)$ (see \cite[Theorem 2.1]{hefez}). %They prove the following result. Let $D_1$ and $D_2$ be two irreducible curves with the same semigroup and the same set of values of K\"ahler differentials. For $i\in\lra{1,2}$, let $\varphi_i$ be a normal parametrization of $D_i$ (see \cite{hefez} for the definition). Then if there exist

\medskip

We first recall some results on equisingular and admissible deformations. We then prove that the stratification by the values of \loga{} residues is finite and constructible, and it refines the stratification by the Tjurina number (see Propositions \ref{finite} and \ref{construc}). We end this section with several algorithms which can be used to compute the set of values of $\mc{R}_D$, inspired by~\cite{bgmescalier} and \cite{hefez-standard}.

\medskip

\begin{de}[\protect{\cite[Definition 10.1.1]{dejong}}]
Let $D$ be a plane curve defined in a neighbourhood $U$ of the origin of $\C^2$ by a reduced equation $f\in\co_{\C^2}(U)$. Let $k\in\N$ and $(S,0)=(\C^k,0)$.  A deformation $F$ of $f$ with base space $S$ is a function $F(x,y,s)\in\co_{\C^2}\widehat{\otimes}\co_S$ which satisfies $F(x,y,0)=f(x,y)$.
\end{de}

For the remainder of this section, we set:

\begin{nota}
Let $X=U\times S$, $\co_X=\co_{\C^2}\widehat{\otimes}\co_S$, $W=F^{-1}(0)\se U\times S$. We assume $F(0,0,s)=0$ for all $s$. For $s\in S$, we set $D_s=W\cap \lrp{ \C^2\times \{s\}}$ and $\m_{S,s}$ the maximal ideal of $\co_{S,s}$, and $F_s=F(.,s)$. In particular, $D_0=D$ and $F_0=f$.
\end{nota}

\subsection{Equisingular and admissible deformations of plane curves}

The following numbers are classical invariants of plane curves:

\begin{de} Let $D$ be a reduced plane curve defined by $f\in\Cxy$.
\begin{itemize}
\item The Milnor number of $f$ is $\mu(f)=\dim_{\C} \Cxy/(f_x',f_y')$
\item The Tjurina number of $f$ is $\tau(f)=\dim_{\C} \Cxy/(f_x',f_y',f)$
\item The delta-invariant of $f$ is $\delta(f)=\dim_\C \co_{\wt{D}}/\co_D$
\end{itemize}
\end{de}

The following proposition gives the relation between $\mu$ and $\delta$:

\begin{prop}[\cite{milnorsing}]
\label{milnor:delta}
We have the following relation: 
$$\mu(f)=2\cdot\delta(f)-p+1,$$
where $p$ is the number of irreducible components of $D$.
\end{prop}

\begin{de}
\label{def:equi}
Let $F$ be a deformation of $f$ with base space $S$. We say that $F$ is an equisingular deformation of $f$ if for all $s\in S$, $\mu(F_s)=\mu(f)$. 
\end{de}

From the equisingularity Theorem for plane curves (see \cite[\S 3.7]{teissier}),  for an equisingular deformation of a plane curve, a parametrization $\varphi$ of $D$ gives rise to a deformation $\varphi_s$ of the parametrization.
We denote by $\val_{D_s}(g)$ the value of $g\in Q(\co_{D_s})$ along $D_s$. Another consequence of the equisingularity Theorem for plane curves is: 

\begin{cor}
\label{cor:equi}
With the same notations, if $F$ is an equisingular deformation of $f$:
\begin{enumerate}
\item All fibers $D_s$ have the same conductor $\gamma$.
\item Let $x(t,s)=(x_1(t_1,s),\ldots,x_p(t_p,s)),y(t,s)=(y_1(t_1,s),\ldots,y_p(t_p,s))$ be a parametrization of~$D_s$.
For all $s\in S$, $$\inf\big(\val_{D_s}(x(t,s)),\val_{D_s}(y(t,s))\big)=\inf\big(\val_D(x(t,0)),\val_D(y(t,0))\big)=(m^{(1)},\ldots,m^{(p)}),$$
where $m^{(j)}$ is the multiplicity of the component $D_j$ of $D$.
\end{enumerate}
\end{cor}
\begin{proof}\ 
\begin{enumerate}
\item By the equisingularity Theorem (see \cite[\S 3.7, (10)]{teissier}), the intersection multiplicity of any couple of branches, and the characteristic exponents of each branch, do not depend on $s$. Thus, the semigroup and therefore the conductor of each branch are also independent from~$s$. Theorem~\ref{cond:delgado} then gives the result. 
%for all $s_1,s_2\in S$, their exists a bijection $b$ between the branches $D_{s_1,1},\ldots,D_{s_1,p}$ of $D_{s_1}$ and the branches $D_{s_2,1},\ldots,D_{s_2,p}$ of $D_{s_2}$ such that for all $i,j\in \unp$, the intersection multiplicity $(D_{s_1,i}.D_{s_1,j})$ is equal to $(D_{s_2,b(i)}.D_{s_2,b(j)})$. In addition, for all $i\in\unp$, the characteristic exponents of $D_{s_1,i}$ and $D_{s_2,b(i)}$ are equal, so that $D_{s_1,i}$ and $D_{s_2,b(i)}$ have the same semigroup, and therefore the same conductor. 
%It comes from Theorem~\ref{cond:delgado}, since by the equisingularity Theorem  (see \cite[\S 3.7, (10)]{teissier}), the conductors and the intersection multiplicities $(D_i.D_j)=\val_i(f_j)$ do not depend on $s$.
\item For all $j\in\unp$, $\inf\big( \val_{D_{s,j}}(x_j(t_j,s)),\val_{D_{s,j}}(y_j(t_j,s))\big)$ is the multiplicity of $D_{s,j}$, which does not depend on $s$ by the equisingularity Theorem. 
\end{enumerate}
\end{proof}

The following proposition will be used in the next subsection, since it gives a common denominator for the \loga{} residues with interesting properties.

\begin{prop}
\label{cst}
There exists $\alpha,\beta\in\C$ such that for all $s$ in a neighbourhood of $0$, $\val(\alpha F_x'(s)+\beta F_y'(s))=\gamma+\lrp{m^{(1)},\ldots,m^{(p)}}-\undun$. In particular, $\alpha F_x'(s)+\beta F_y'(s)$ induces a non zero divisor in $\co_D$ whose value does not depend on $s$. 
\end{prop}
\begin{proof}
Thanks to the equisingularity Theorem, one can prove that  there exists a linear change of coordinates $(u,v)$ such that for all $s$ in a neighbourhood of $0\in S$, $\val_{D_s}(u)=\lrp{m^{(1)},\ldots,m^{(p)}}$. The conclusion follows from Corollary~\ref{cor:equi} and Lemma~\ref{val:jc:xy}.
\end{proof}

We want now to understand the behaviour of a generating family of the module of residues.

\medskip

We recall that plane curves are free divisors. Moreover, they are the only free divisors with isolated singularities, since by \cite{alek}, the singular locus of a free divisor is of codimension one in the hypersurface. The equisingularity assumption is not sufficient to obtain a deformation $(\rho_1(s),\rho_2(s))$ of a generating family of $\mc{R}_D$ such that $(\rho_1(s),\rho_2(s))$ generate $\mc{R}_{D_s}$: equisingularity is not the "good" functor of deformation for free divisors. A functor of deformation adapted to free divisors is suggested by M.Torielli in \cite{tor}.

The following definition is equivalent to the definition of M.Torielli (see \cite[Definition 3.1]{tor}) thanks to both \cite[Proposition 3.7]{tor}
and \cite[Theorem 1.91]{greueldefo}:
\begin{de}
\label{def:adm}
Let $D$ be a free divisor defined in a neighbourhood of $0\in\C^n$ by a reduced equation $f$. An admissible deformation $X$ of $D$ with base space $S$ is a deformation of $D$  such that the module $\co_{\C^n\times S,0}/(F,F_x',F_y')$ is a flat $\co_{S,0}$-module.
\end{de}

The following proposition describes an admissible deformation of a plane curve thanks to the Tjurina number.

\begin{prop}
Let $F$ be a deformation of $f\in\co_{\C^2}(U)$ with base space $S$ such that for all $s\in S$, $\sum_{x_j\in Sing(D_s)} \tau_{x_j}=\tau_0$ where $\tau_0$ is the Tjurina number of $D_0$. Such a deformation is an admissible deformation.
\end{prop}
\begin{proof}
 We denote by $p$ the restriction of the canonical epimorphism $\C^2\times S\to S$ to the relative singular locus, which is a finite morphism. We set $\mc{F}=p_*\lrp{\co_{\C^2\times S}/\lrp{F,\frac{\dr F}{\dr x},\frac{\dr F}{\dr y}}}$. Then $\mc{F}_s/ \m_{S,s} \mc{F}_s= \co_{\C^2}/\lrp{F_s, \frac{\dr F_s}{\dr x},\frac{\dr F_s}{\dr y}}$, which is by assumption a complex vector space of dimension $\tau_0$. The proposition is then a direct consequence of \cite[Theorem 1.81]{greueldefo}.
\end{proof}
%We use Theorem $1.81$ of~\cite{greueldefo} for $\mc{F}=p_*\lrp{\co_{\C^2\times S}/(F,F_x',F_y')}$ with $p : \C^2\times S\to S$ the canonical epimorphism (which is finite). The map $s\mapsto \dim_\C \mc{F}_s/\m_{S,s}\mc{F}_s$ is then constant, thus, $\mc{F}_0$ is a flat $\co_{S,0}$-module. Since $\mc{F}_0=\co_{\C^2\times S,0}/(F,F_x',F_y')$, the result follows.% and $p : \times S\to S$ the canonical epimorphism.  

\begin{prop}[\protect{\cite[Lemma 3.22]{tor}}]
\label{chp:def}
Let $F(x,y,s)$  be an equisingular and admissible deformation of the plane curve defined by $f$ with base space $S$. Let $(\delta_1,\delta_2)$ be a basis of the module of \loga{} vector fields along $D$. Then $\delta_1,\delta_2$ induce relations between $f,f_x',f_y'$. By flatness, we can extend them to obtain relative \loga{} vector fields $\wt{\delta}_1,\wt{\delta}_2\in (\Theta_{U\times S/S})/(\m_S \Theta_{U\times S/S})$ of $F$. Then, for $s$ in a neighbourhood of $0\in S$, $\lrp{\wt{\delta}_1(s),\wt{\delta}_2(s)}$ is a basis of $\der(-\log D_s)$. 
\end{prop}

\begin{cor}
\label{res:gen}
Let $\wt{\delta}_i=A_i(s)\dr_x+B_i(s)\dr_y$, $i=1,2$ be as in Proposition~\ref{chp:def}. Thanks to the duality between the modules $\mathrm{Der}(-\log D_s)$ and $\Omega^1(\log D_s)$, the following elements  generate the module of residues for all $s$ in a neighbourhood of $0\in S$:

$$\left\{
\begin{aligned}
\rho_1(s)&=\f{-\beta A_2(s)+\alpha B_2(s)}{\alpha F_x'(s)+\beta F_y'(s)}\\
\rho_2(s)&=\f{\beta A_1(s)-\alpha B_1(s)}{\alpha F_x'(s)+\beta F_y'(s)}
\end{aligned}
\right.$$ 
where $\alpha,\beta\in\C$ are given by Proposition~\ref{cst}.
\end{cor}

\subsection{Properties of the stratification by \loga{} residues}
We consider an equisingular deformation $F$ of $f$ with base space  $(S,0)\simeq (\C^k,0)$ for a $k\in\N$. We denote by $\mc{R}_s$ the module of \loga{} residues of $D_s$. 

The purpose of this section is to study the stratification by the values of logarithmic residues defined in Definition \ref{strat}. We prove that this stratification is finite, finer than the stratification by the Tjurina number and constructible (see Propositions \ref{finite} and \ref{construc}). We end with two examples, the first one shows that the stratification by logarithmic residues may be strictly finer than the stratification by the Tjurina number (see Example \ref{res:tau}). The stratification by logarithmic residues of Example \ref{front} does not satisfy the frontier condition.

\begin{de}
\label{strat}
Let $F(x,y,s)$ be an equisingular deformation of a reduced plane curve $D$ with $p$ branches defined by $f\in\Cxy$, with base space $S$. The stratification by \loga{} residues is the partition $S=\bigcup_{\mc{V}\subseteq \Z^p} S_{\mc{V}}$  where $s\in S_{\mc{V}}$ if and only if $\val(\mc{R}_s)=\mc{V}$. 
\end{de}
An example of a stratification by the values of logarithmic residues is given in example~\ref{res:tau}.

\smallskip

Let us compare the stratification by logarithmic residues with the stratification by the Tjurina number. The stratification by the Tjurina number is the partition $S=\bigcup_{n\in\N} S_n$ where $s\in S_n$ if and only if $\tau(F_s)=n$. This stratification is finite since the Tjurina number is bounded by the Milnor number, which is constant on $S$ by the equisingularity condition.

\begin{prop}
\label{prop:dim:res}
Let $D$ be a plane curve germ. Then:
\begin{equation}
\label{dim:res}
\dim_\C \mc{R}_D/\co_{\wt{D}}=\tau-\delta.
\end{equation}
\end{prop}

\begin{proof}
Thanks to Propositions~\ref{duality:MCM} and~\ref{th:dual:jd} we have:
\begin{align*}
\dim_\C \mc{R}_D/\co_{\wt{D}} &=\dim_\C \mc{R}_D/\co_D-\dim_\C \co_{\wt{D}}/\co_D=\dim_\C \co_D^\vee/\mc{R}_D^\vee-\delta=\dim_\C \co_D/\mc{J}_D-\delta\\
 &=\tau-\delta
\end{align*}
\end{proof}

\begin{prop} 
\label{finite}The stratification by \loga{} residues satisfies the following properties:
\begin{enumerate}
\item The stratification by logarithmic residues is finer  than the stratification by the Tjurina number.
\item The stratification by logarithmic residues is finite. 
\end{enumerate}
\end{prop}
\begin{proof}
The first claim is a direct consequence of Proposition~\ref{prop:dim:res}, since the equisingularity condition ensures that $\delta(F_s)$ does not depend on $s$, and the dimension of the quotient $\mc{R}_s/\co_{\wt{D_s}}$ can be computed from the values of $\mc{R}_s$ by Corollary~\ref{suite}. The second claim comes from both Proposition~\ref{cst}, which gives a lower bound $u$ of the set of values of \loga{} residues which do not depend on~$s$, and Corollary~\ref{neg:res}. As a consequence, the values of $\mc{R}_s$ are determined by the values $v$ of $\mc{R}_s$ satisfying $u\leqslant v\leqslant 0$. 
\end{proof}

\begin{prop}
\label{construc}
Each stratum $S_{\mc{V}}$ of the stratification by \loga{} residues is constructible. If moreover $D$ is irreducible, then each stratum is locally closed. 
\end{prop}

The hypothesis of $D$ being irreducible was forgotten in  the corresponding statement \cite[Proposition 4.2]{polcras}.
\begin{proof}
%For lack of reference, we suggest a proof.  
 By the appendix by Teissier in~\cite{zariskicourbe}, the strata of the stratification by the Tjurina number are locally analytic and locally closed. It is therefore sufficient to consider the behaviour of logarithmic residues in a $\tau$-constant stratum $S_\tau$. For the sake of simplicity, we denote $S=S_\tau$.  
 
 By Corollary~\ref{res:gen}, for all $s$, the  $\co_S$-module $\mc{R}_s$ is generated by 
 $$\left\{
\begin{aligned}
\rho_1(s)&=\f{-\beta A_2(s)+\alpha B_2(s)}{ \alpha F_x'(s)+\beta F_y'(s)}\\
\rho_2(s)&=\f{\beta A_1(s)-\alpha B_1(s)}{\alpha F_x'(s)+\beta F_y'(s)}
\end{aligned}
\right.$$ 
where $\alpha,\beta\in\C$ are given by Proposition~\ref{cst}. The value of the common denominator $\alpha F_x'(s)+\beta F_y'(s)$ does not depend on $s$, so that it is sufficient to consider the values of the numerators.

\medskip

We denote by $N_1$ and $N_2$ the numerators of $\rho_1(s)$ and $\rho_2(s)$.   We recall that the values $v$ of $\mc{R}_s$ satisfying $v\leqslant 0$ are sufficient to determine $\val(\mc{R}_s)$, so that it is sufficient to consider the  set of all elements $x^iy^jN_k$ with $i,j\in\N, k\in\lra{1,2}$ such that $\val(x^iy^jN_k)\leqslant u$, where $u=\val(\alpha F_x'(s)+\beta F_y'(s))$. We set $\lra{X_1,\ldots,X_q}:=\lra{x^iy^jN_k; \val(x^iy^jN_k)\leqslant u}$, where the elements are numbered in an arbitrary order, and $q$ is the number of elements in this set. 

For all $i\in\lra{1,\ldots,q}$, we have $X_i=\lrp{\sum_{j\geqslant 0} a_{i,j,1}(s) t_1^j,\ldots,\sum_{j\geqslant 0} a_{i,j,p}(s) t_p^j}\in \co_{\wt{D_s}}$. 

For $v\in\Z^p$ and $k\in\unp$ we set $X_{i,k}^{v}(s)=(a_{i,0,k}(s),\ldots, a_{i,v_k,k}(s))\in \co_S^{v_k+1}$. For $v\in\Z^p$ we define the following matrix $A_v(s)\in\mc{M}_{q,\ell_v}(\co_S)$ with $\ell_v=\sum_{j=1}^p (v_j+1)$:

$$A_v(s)=\begin{pmatrix}
(X_{1,1}^{v}(s)) & \ldots & (X_{1,p}^v(s)) \\
\vdots & & \vdots \\
(X_{q,1}^v(s)) &\ldots & (X_{q,p}^v(s))
\end{pmatrix}.$$

The $i$th row of this matrix encodes the respective Taylor developments of $X_i$ along the branch $D_k$ for $k=1,\cdots, p$ up to order $v_k$.

We set $(e_1,\ldots, e_p)$ the canonical basis of $\Z^p$. 
We use the rank of the matrices $A_v(s)$ to characterize the property $v\in\val(\mc{R}_s)$ for $s\in S$:

$$v\in\val(\mc{R}_{s})\iff \forall k\in\lra{1,\ldots,p}, \mathrm{rank}\lrp{A_{v-\undun}(s)}< \mathrm{rank}\lrp{A_{v-\undun +e_k}(s)}.$$

Indeed, if the conditions of the right-hand  side are satisfied, then for all $k\in\lra{1,\ldots,p}$, there exists a linear combination $M_k=\sum_{i=1}^q \lambda_{i,k} X_i(s)$ with $\lambda_{i,k}\in\C$ such that $\val(M_k)\geqslant v$ and $\val_k(M_k)=v_k$. We use Proposition \ref{propinf} to conclude. 

Therefore, for a given $\mc{V}\subseteq \Z^p$ for which $S_{\mc{V}}\cap S \neq \emptyset$ and $\ov{\mc{V}}:=\lrp{\mc{V}+u}\cap\lra{w\in\Z^p; 0\leqslant w\leqslant u}$:

$$s\in S_{\mc{V}}\iff s\in \bigcap_{v\in\ov{\mc{V}}} \lrp{\bigcup_{1\leqslant r\leqslant M} \lrp{V\lrp{\mc{F}_r\lrp{A_{v-\undun}(s)}}\cap \bigcap_{1\leqslant k\leqslant p} \lrp{V\lrp{\mc{F}_r(A_{v-\undun+e_k}(s)}^c}}},$$
where $\mc{F}_r(A)$ denotes the ideal generated by the $r\times r$ minors of the matrix $A$ and $M=\min(q,\ell_v+1)$. We notice that the elements $v\notin \mc{V}$ can not be reached since otherwise, by Corollary~\ref{valeq}, the dimension of $\mc{R}_{s}/\co_{D_s}$ would be strictly greater than $\tau-\delta$.%, which can be proved by an argument similar to the first two steps of the proof of Theorem~\ref{symmetry-values}.

Hence the result for reducible curves.

\smallskip

Let us assume now that $D$ is irreducible. In this case, the rank of the matrix $A_v(s)$ increases exactly by $1$ when a valuation is reached. We set $\ov{\mc{V}}=\lra{v_1<\ldots< v_L}=\lrp{\mc{V}+u}\cap\lra{0,\ldots,u}$. 
Then:

$$s\in S_{\mc{V}}\iff s\in \bigcap_{\ell=1}^L \bigcap_{j=v_{\ell-1}+1}^{v_\ell-1} \big(V(\mc{F}_{\ell}(A_{j}(s))\cap\lrp{V(\mc{F}_{\ell}(A_{v_\ell}}^c\big).$$

Therefore the stratum $S_{\mc{V}}$ is locally closed. 
\end{proof}

We recall here the examples of \cite{polcras} with more details. The first example shows that the stratification by \loga{} residues may be strictly finer than the stratification by the Tjurina number, whereas the second shows that the stratification by \loga{} residues does not satisfy the frontier condition defined below.

\begin{ex}
\label{res:tau}
We consider $f(x,y)=x^5-y^6$ and the equisingular deformation of $f$ given by $F(x,y,s_1,s_2,s_3)= x^5-y^6+s_1x^2y^4+s_2x^3y^3+s_3x^3y^4$. The stratification by $\tau$ is composed of three strata, $S_{20}=\{0\}$, $S_{19}=\{(0,0,s_3),s_3\neq 0\}$ and $S_{18}=\{(s_1,s_2,s_3), (s_1,s_2)\neq (0,0)\}$, where the index indicates the value of $\tau$. The computation of the values of $\mc{J}_{D_s}$ is quite easy in this case and gives thanks to Theorem~\ref{symmetry-values} the values of $\mc{R}_{D_s}$. The stratification of $\C^3$ by the values of $\mc{R}_{D_s}$ is then $\C^3=S_{20}\sqcup S_{19}\sqcup S_{18}'\sqcup S_{18}''$  where $S_{18}'=\lra{(s_1,s_2,s_3),s_1\neq 0}$ and $S_{18}''=\lra{(0,s_2,s_3),s_2\neq0}$. The stratum $S_{18}$ splits up into the two strata  $S_{18}'$ and $S_{18}''$ for stratification by the values of $\mc{R}_{D_s}$, and the other strata of the stratification by the Tjurina number are the same as by the values of logarithmic residues. Therefore, the stratification by logarithmic residues is strictly finer than the stratification by $\tau$. The corresponding values are:

\medskip

\begin{tabular}{|c|c|l|}
\hline
Stratum & $\dim_\C \mc{R}_{D_s}/\co_{\wt{D_s}}$ &  negative values of $\mc{R}_{D_s}$\\
\hline
$S_{20}$ & 10 & $-1,-2,-3,-4,\ \ -7,-8,-9,\ \ \ -13,-14,\ \ \ \ -19$\\
\hline
$S_{19}$ & 9 & $-1,-2,-3,-4,\ \ -7,-8,-9,\ \ \ -13,-14\ \ \ \ \ $ \\
\hline
$S_{18}'$ & 8 & $-1,-2,-3,-4,\ \ -7,-8,-9,\ \ \ \ \ \ \ \ \ \ -14\ \ \ \ \ $\\
\hline
$S_{18}''$ & 8 & $-1,-2,-3,-4,\ \ -7,-8,-9,\ \ \ -13\ \ \ \ \ \ $\\
\hline
\end{tabular}

 \medskip

\end{ex}

Example~\ref{front} below shows that the stratification by logarithmic residues do not necessarily satisfy the frontier condition defined below. 
\begin{de}
A stratification $S=\bigcup_{\alpha} S_\alpha$ satisfies the \emph{frontier condition} if for $\alpha\neq \beta$, $S_\alpha\cap \overline{S_\beta}\neq\emptyset$ implies $S_\alpha\subseteq \overline{S_\beta}$, with $\overline{S_\beta}$ the closure of $S_\beta$.
\end{de}

We first prove the following property:

%  A stratification $S=\bigcup_{\alpha} S_\alpha$ satisfies the frontier condition if for $\alpha\neq \beta$, $S_\alpha\cap \overline{S_\beta}\neq\emptyset$ implies $S_\alpha\subseteq \overline{S_\beta}$, with $\overline{S_\beta}$ the closure of $S_\beta$.

\begin{lem}
Let $D$ be a quasi-homogeneous plane curve germ with conductor $\gamma$. Then: $$\gamma-\undun+\lrp{\val(\co_D)\backslash \{0\}}=\val(\mc{J}_D).$$
\end{lem}
  
\begin{proof}
Let $p$ be the number of branches.
The inclusion $\subseteq$ is given by Remark~\ref{jac:semi}. For the other inclusion, we notice that  Remark~\ref{jac:semi} implies $t^{2\gamma-\undun}\co_{\wt{D}}\subseteq\mc{J}_D\subseteq \mc{C}_D$. We have the following equality:
$$\dim_\C \mc{C}_D/t^{2\gamma-1}\co_{\wt{D}}=\dim_\C \mc{C}_D/\mc{J}_D+\dim_\C \mc{J}_D/t^{2\gamma-\undun}\co_{\wt{D}}.$$

By Propositions \ref{prop:dim:res} and \ref{duality:MCM}, we have $\dim_\C  \mc{C}_D/\mc{J}_D=\tau-\delta$. Since $D$ is quasi-homogeneous, we have $\tau=\mu$ so that by Proposition \ref{milnor:delta} we have $\dim_\C \mc{C}_D/\mc{J}_D=\delta-p+1$. Moreover, since $\delta=\dim_\C \co_{\wt{D}}/\co_D$, and since $C$ is Gorenstein, $\dim_\C \co_D/\mc{C}_D=\delta$, we have $\dim_\C \co_{\wt{D}}/\mc{C}_D=2\delta$, thus $\dim_\C \mc{C}_D/t^{2\gamma-1}\co_{\wt{D}}=2\delta-p$. Therefore:
$$\dim_\C \mc{J}_D/t^{2\gamma-1}\co_{\wt{D}}=\delta-1.$$
 Let $\m$ be the maximal ideal of $\co_D$. Then $\val(\m)=\val(\co_D)\backslash \{0\}$ and the quotient $\m/\mc{C}_D$ has dimension $\delta-1$. Therefore, $\dim_\C t^{\gamma-\undun} \m/t^{2\gamma-\undun}\co_{\wt{D}}=\delta-1$. Since $\val(t^{\gamma-\undun} \m)\subseteq \val(\mc{J}_D)$, by Corollary~\ref{valeq}, the equality follows. 
\end{proof}  
  % If $\val(\mc{J}_D)\neq \gamma-\undun+\lrp{\val(\co_D)\backslash \{0\}}$, by an argument similar to the one of the two first steps of the proof of Theorem~\ref{symmetry-values}, we have $\dim_\C \mc{J}_D/t^{2\gamma-1}\co_{\wt{D}}>\delta-1$, which is a contradiction. Hence the result.
  \begin{ex}
\label{front}
Let us consider the deformation $F(x,y,s_1,s_2)=x^{10}+y^8+s_1x^5y^4+s_2x^3y^6$ for $s_1,s_2$ in a neighbourhood of $0$ so that the deformation is equisingular. It is given in~\cite{bmgbfundefo}, as an example of the stratification by the $b$-function not satisfying the frontier condition.

 Contrary to the previous example, this curve is not irreducible. %We first give a property of quasi-homogeneous curves:
%Let us come back to our example.

We  notice that $F(x,y,s_1,0)$ is quasi-homogeneous. Therefore, the previous lemma shows that the values of the Jacobian ideal along the quasi-homogeneous stratum does not change. Therefore, the quasi-homogeneous stratum is a stratum of the stratification by \loga{} residues. 

Moreover, one can check that there are three strata for the stratification by the Tjurina number: the  quasi-homogeneous stratum $S_1$ defined by $s_2=0$ for which $\tau=63$, a stratum $S_2$ defined by $s_1=0$ and $s_2\neq 0$ for which $\tau=54$ and a stratum $S_3$ defined by $s_1s_2\neq 0$ for which $\tau=53$. Therefore, the stratification by \loga{} residues does not satisfy the frontier condition. Indeed, the stratification by logarithmic residues is finite and constructible so that there is a stratum of the stratification by logarithmic residues which is an open dense subset of $S_2$, which therefore contains the origin in its closure, but not the whole quasi-homogeneous stratum.
\end{ex}

\subsection{Algorithms to compute the \loga{} residues along plane curves with one or two components}

We suggest here several methods which can be used to compute the values of \loga{} residues. 

Thanks to the symmetry Theorem~\ref{symmetry-values}, computing the values of $\mc{J}_D$ is equivalent to the computation of the values of $\mc{R}_D$.

\subsubsection{Irreducible semi-quasi homogeneous polynomials}

This algorithm is used to study the equisingular deformation of a quasi-homogeneous polynomial of the form $x^a-y^b$, with $\gcd(a,b)=1$. It is inspired by \cite{bgmescalier}.

\s 

We consider the following equation of an irreducible curve, with $s_{ij}\in\C$ and $\gcd(a,b)=1$:
\begin{equation}
\label{defo:ab}
F(x,y)=x^a-y^b +\sum_{\substack{1\leqslant i<a-1 \\ 1\leqslant j<b-1 \\ ib+ja>ab}} s_{ij}x^iy^j.
\end{equation}
A parametrization of the curve is given by $x(t)=t^b+g(t), y(t)=t^a+ h(t)$ where $g,h\in\C\{t\}$ with $\val(g)>b, \val(h)>a$. 

\smallskip

We set for $i,j\in\N^2$, $\rho(i,j)=ib+ja$. We define a monomial ordering by: $(i,j)<(i',j')$ \ssi{} $\rho(i,j)<\rho(i',j')$ or $\big(\rho(i,j)=\rho(i',j') \text{ and } i<i'\big)$. If $H=\sum_{i,j} a_{i,j} x^iy^j\in\Cxy$ is non zero, we set $\exp(H)=\min\lrp{(i,j), a_{i,j}\neq 0}$ and $\rho(H):=\rho(\exp(H))$.

Polynomials of the form~\eqref{defo:ab} are studied in \cite{bgmescalier}. The authors give an algorithm to compute the "\textit{escalier}" of the curve, which is by definition the complement in $\N^2$ of the set $$E=\{\exp(g);g\in(F,F_x',F_y')\se\Cxy\}.$$

More precisely, they give the explicit computation of a finite family $(A_j)_{-1\leqslant j\leqslant K}$ of points of $\N^2$ such that $E=\bigcup_{j=-1}^{K} A_j+\N^2$, and none of the $A_j$'s can be removed.  Then it is possible to prove: 

\begin{prop}
We have the following equality: $$\val(\mc{J}_D)=\bigcup_{i=-1}^{K} \Big(\rho(A_j)+\val(\co_D)\Big).$$
\end{prop}

\subsubsection{Irreducible plane curve}
\label{irred:algo}

In \cite{hefez-standard}, an algorithm is proposed to compute the set of values of K\"ahler differentials of an irreducible plane curve.  By Proposition~\ref{res:kahler}, it gives also the values of $\mc{R}_D$. In fact, one can see that the algorithm of \cite{bgmescalier} corresponds to the algorithm of \cite{hefez-standard} by Proposition~\ref{jc:kahler}.  

 Moreover, if a generating family of $\mc{R}_D$ is known, the algorithm of \cite[Theorem 2.4]{hefez-standard} can be used directly on this family to compute a standard basis $(H,G)$ for $\mc{R}_D$ (see \cite[Definition 2.1]{hefez-standard}). In particular, $G$ is a standard basis of $\co_D$ and $H\subseteq \mc{R}_D$.  In order to determine $\val(\mc{R}_D)$ from $H$ and $G$, we need the following notion: 
 
 \begin{de}[\protect{\cite{hefez-standard}}]
 A $G$-product is an element of the form $\prod_{i=1}^q g_i^{\alpha_i}$ with $q\in\N$, $\alpha_j\in\N$ and $g_i\in G$.
 \end{de}

A standard basis $G$ of $\co_D$ is characterized by the fact that for all  $h\in\co_D$, there exists a $G$-product $g$ such that $\val(h)=\val(g)$ (see \cite{hefez-standard}).
\begin{ex}
Let $D$ be the  irreducible plane curve parametrized by $x(t)=t^4$ and $y(t)=t^6+t^7$. By \cite[Example 5.2.13]{dejong}, the semigroup of $D$ is generated by $4,6,13$ so that a  standard basis of $\co_D$ is  $G=\lra{x,y,y^2-x^3}$. The $G$-products are then the elements $x^iy^j(y^2-x^3)^k$ for $i,j,k\in\N$.  
\end{ex}

For any irreducible curve $D$, if $(H,G)$ is a standard basis of $\mc{R}_D$, then: 
$$\val(\mc{R}_D)=\lra{\val(h)+\val(g); h\in H, g \text{ a } G\text{-product }}.$$

\subsubsection{Plane curves with two branches}

Let $D=D_1\cup D_2$ be a plane curve germ with two irreducible components.  We suggest here an algorithm to compute the set of negative values of $\mc{R}_D$. It is more technical than in the irreducible case, and cannot be generalized to plane curves with three or more branches. It can be compared to the fact that the analytic classification proposed in \cite{Hefez2} for two branches is also more complicated than in the irreducible case, and can not  be easily extended to plane curves with three or more branches.

Example~\ref{ex:algo} illustrates the algorithm for two branches which is suggested below.

\smallskip

%The algorithm in \cite{hefez-standard} is given for irreducible curves, for which the set of valuations is totally ordered, so that we cannot apply it directly to reducible plane curves. Nevertheless, we can use it if we consider only one of the components.  We denote for $p\in\lra{1,2}$, $\Z_{\leqslant 0}^p=\lra{v\in\Z^p; \ v\leqslant 0}$.

\begin{remar}
\label{algo}
The algorithm in \cite{hefez-standard} is given for irreducible curves, for which the set of valuations is totally ordered, so that we cannot apply it directly to reducible plane curves. Nevertheless, we can use it if we consider only one of the components.  More precisely, if we consider an ideal $I=(h_1,\ldots,h_q)$ in $\co_D$, we associate to it a $\C\lra{x_1(t_1),y_1(t_1)}$-module $\ov{I}=\lrp{\ov{h}_1,\ldots,\ov{h}_q}\subset \C\lra{t_1}$, where $\ov{h}_i$ is the image of $h_i$ in $\C\lra{t_i}$. The algorithm of \cite[Theorem~2.4]{hefez-standard} applied to $\ov{I}$ gives a standard basis $(\ov{H}_1,\ov{G}_1)$ for $\ov{I}$. This algorithm is based on the notion of $S$-process (see \cite[Definition~2.2]{hefez-standard}), so that we can simultaneously compute a family $(H_1, G_1)$ in $\co_D$ such that the image of $G_1$ in $\C\lra{t_1}$ is $\ov{G}_1$, and the image of $H_1$ in $\C\lra{t_1}$ is $\ov{H}_1$. 
\end{remar}

We denote for $p\in\lra{1,2}$, $\Z_{\leqslant 0}^p=\lra{v\in\Z^p; \ v\leqslant 0}$.

\s

\etape{First step}

\s

First of all, we set $g\in (f_x',f_y')$ a non zero divisor of $\co_D$, and we fix it as the common denominator of all residues of $D_1$, $D_2$ and $D$, so that we can consider only the numerators to compute the set of values in each case. 

\smallskip

Let $i\in\lra{1,2}$. We consider only the branch $D_i$. Thanks to the algorithm of \cite[Theorem~2.4]{hefez-standard} and Remark~\ref{algo}, we compute (thanks to the numerators) a family $(R_i, G_i)$  in $Q(\co_D)$ such that its image in $\C\lra{t_i}$ is a standard basis of $\restriction{\mc{R}_{D_i}}{D_i}$.

To determine entirely $\val(\mc{R}_D)\cap \Z_{\leqslant 0}^2$, we first compute the projection $\val_1(\mc{R}_D)$ of $\val(\mc{R}_D)$. To do this, we apply again \cite[Theorem~2.4]{hefez-standard} and Remark~\ref{algo}  to obtain a family $(R, G_1)$ in $Q(\co_D)$ such that its image in $\C\lra{t_1}$ is a standard basis of $\restriction{\mc{R}_D}{D_1}$.  Therefore, for all  $v_1\in\val_1(\mc{R}_D)\cap (\val_1(\mc{R}_{D_1}))^c$, there exists $\rho\in R$ such that $\val_1(\rho)=v_1$. 

%Let $i\in\lra{1,2}$. We consider only the branch $D_i$. By applying the algorithm of \cite[Theorem~2.4]{hefez-standard} we compute a standard basis $G_i$ of $\co_{D_i}$. We use it to compute also a standard basis $R_i$ of $\mc{R}_{D_i}$, thanks to a generating family of $\mc{R}_{D_i}$. %As in \cite{hefez-standard}, a $G_i$-product is an element of the form $\prod_{j=1}^q h_j^{\alpha_j}$ where $h_j\in G_i$, $\alpha_j\in\N$ and $q\in\N$.

\s

%To determine entirely $\val(\mc{R}_D)\cap \Z_{\leqslant 0}^2$, we first compute the projection $\val_1(\mc{R}_D)$ of $\val(\mc{R}_D)$. To do this, we apply the algorithm of \cite[Theorem 2.4]{hefez-standard} to a generating family of $\mc{R}_D$, but by considering only the valuation along $D_1$.  We deduce a family $R$ of elements of $\mc{R}_D$ such that for all $v_1\in\val_1(\mc{R}_D)\cap (\val_1(\mc{R}_{D_1}))^c$, there exists $\rho\in R$ such that $\val_1(\rho)=v_1$. 

%To determine entirely $\val(\mc{R}_D)\cap \Z_{\leqslant 0}^2$, we first compute the projection $\val_1(\mc{R}_D)$ of $\val(\mc{R}_D)$. To do this, we consider the image in $\C\lra{t_1}$ of a generating family of $\mc{R}_D$. Thanks to the algorithm of \cite[Theorem~2.4]{hefez-standard}, we obtain a  family $R$ of elements of $\mc{R}_D$ such that for all $v_1\in\val_1(\mc{R}_D)\cap (\val_1(\mc{R}_{D_1}))^c$, there exists $\rho\in R$ such that $\val_1(\rho)=v_1$. %and we apply the algorithm of \cite[Theorem 2.4]{hefez-standard}. to a generating family of $\mc{R}_D$, but by considering its image in $\C\lra{t_1}$.   We deduce a family $R$ of elements of $\mc{R}_D$ such that for all $v_1\in\val_1(\mc{R}_D)\cap (\val_1(\mc{R}_{D_1}))^c$, there exists $\rho\in R$ such that $\val_1(\rho)=v_1$. 

\s

\etape{Second step}

\s

We set $\mc{M}_0=\lra{(0,v_2); v_2\in\val_2(\mc{R}_{D_2})\cap\Z_{\leqslant 0}}$ and $H_0=R_2$. By Proposition~\ref{zerodiv:res}, it gives all the values of $\val(\mc{R}_D)\cap \big(\lra{0}\times\Z_{\leqslant 0}\big)$. 

\smallskip

Let us assume that for a $k\in\N^*$ we have constructed  sets $\mc{M}_{k-1}$ and $H_{k-1}\subseteq \mc{R}_D$ such that $$\mc{M}_{k-1}=\lra{(v_1,v_2)\in\val(\mc{R}_D); -k+1\leqslant v_1\leqslant 0 \text{ and } v_2\leqslant 0},$$ and for all $v_2\in\val_2(\mc{M}_{k-1})$, there exists $\rho\in H_{k-1}$ and a $G_2$-product $h$ with $\val_2(h\cdot\rho)=v_2$ and $\val_1(h\cdot\rho)\geqslant -k+1$. 

\smallskip

Let us compute $\mc{M}_{k}$ and $H_{k}$. If $-k\notin\val_1(\mc{R}_D)$, $\mc{M}_{k}=\mc{M}_{k-1}$ and $H_{k}=H_{k-1}$. Otherwise, there are several cases to consider. 

\smallskip

\textit{First case: $-k\in\val_1(\mc{R}_{D})\cap\val_1(\mc{R}_{D_1})$}. It means that $(-k,\infty)\in\ov{\val(\mc{R}_D)}$, which is by Proposition~\ref{zero:div} equivalent to $ (-k,0)\in \val(\mc{R}_D)$. By Proposition~\ref{propinf}, one can see that \begin{equation}
\label{mk}
\mc{M}_{k}\supseteq\mc{M}_{k-1}\cup \lra{(-k,v_2), v_2\in\val_2(\mc{M}_{k-1})}.
\end{equation} Moreover, by Proposition~\ref{valquimonte}, if $(-k,v_2)\in\val(\mc{R}_D)$ with $v_2\leqslant 0$, then $v_2\in\val_2(\mc{M}_{k-1})$. Therefore the inclusion in \eqref{mk} is an equality and $H_{k}=H_{k-1}$.

\smallskip

\textit{Second case:} $-k\in\val_1(\mc{R}_D)$ but $-k\notin \val_1(\mc{R}_{D_1})$. There exists $\rho_0\in R$ such that $\val_1(\rho_0)=-k$. Let $w_2=\val_2(\rho_0)$. We may assume by Proposition~\ref{propinf} that $w_2\leqslant 0$ since $0\in\val( \mc{R}_D)$. 

\smallskip

%Thanks to Propositions \ref{propinf} and \ref{valquimonte}, one can check that: 

\indent \indent \textit{First sub-case: $w_2\notin \val_2(\mc{M}_{k-1})$}. Then:
\begin{equation}
\label{M}
\mc{M}_{k}=\mc{M}_{k-1}\cup\lra{(-k,w_2)}\cup\lra{(-k,v_2); v_2\in\val_2(\mc{M}_{k-1}) \text{ and } v_2\leqslant w_2}
\end{equation}
and $H_{k}=H_{k-1}\cup\lra{\rho_0}$. Indeed, the inclusion $\supseteq$ of~\eqref{M} comes from Proposition~\ref{propinf}. By Proposition~\ref{valquimonte} and an argument similar to the argument of the first case, one can prove the equality in~\eqref{M}. 

\smallskip

\indent \indent \textit{Second sub-case: $w_2\in\val_2(\mc{M}_{k-1})$}. Thanks to Propositions \ref{propinf} and \ref{valquimonte}, one can check that by a convenient linear combination of $\rho_0$ and elements of form $h\cdot \rho$ with $h$ a $G_2$-product and $\rho\in H_{k-1}$, there exists $\rho_0'\in\mc{R}_D$ with $\val(\rho_0')=(-k,w_2')$ and $w_2'\notin \val_2(\mc{M}_{k-1})$. We then recognize the previous sub-case, and we have $H_{k}=H_{k-1}\cup\lra{\rho_0'}$.

We can stop when the minimal value $-q$ of $\val_1(\mc{R}_D)$ is reached. Then, by Proposition \ref{neg:res}:
$$\val(\mc{R}_D)=\mc{M}_q\cup\lra{v\in\Z^p; \inf(v,0)\in\mc{M}_q}.$$

\begin{ex}
\label{ex:algo}
Let us consider the plane curve $D$ defined by $f=f_1f_2$ with $f_1=y^3-x^5$ and $f_2=y^5-x^3$. We denote by $D_1$ the curve defined by $f_1$, and $D_2$ the curve defined by $f_2$. A parametrization of the curve $D$ is $x=(t_1^3,t_2^5), y=(t_1^5,t_2^3)$.  Computations show that  generators of the module of logarithmic residues $\mc{R}_D$ are (see \cite[\S 6.2.4]{polthese}):
$\rho_1=-\frac{15(x^{3}y^{4}-xy^{2})}{f_y'}$ and $\rho_2=\frac{-9x^3+25x^5y^{2}-16y^5}{f_y'}$.  As elements of $Q(\co_D)=Q(\co_{D_1})\oplus Q(\co_{D_2})$, they give:
$\rho_1=\lrp{\frac{-5}{t_1^6},\frac{3}{t_2^{10}}}$ and $\rho_2=\lrp{\frac{3}{t_1^{10}},\frac{-5}{t_2^6}}$.

Since the branches are quasi-homogeneous, one can see with Saito criterion (see \cite{saitolog}) that $\frac{\dd f_1}{f_1}$ and $\frac{5y\dd x-3x\dd y}{f_1}$ is a basis of $\Omega^1(\log D_1)$, and $\frac{\dd f_2}{f_2}$ and $\frac{3y\dd x-5x \dd y}{f_2}$ is a basis of $\Omega^1(\log D_2)$. Therefore, the generators of $\mc{R}_{D_1}$ are $\rho_{11}=1$ and $\rho_{12}=\frac{-3x}{{\dr_y f_1}}$, where $\dr_y f_1=\frac{\dr f_1}{\dr y}$, and the generators of $\mc{R}_{D_2}$ are $\rho_{21}=1$ and $\rho_{22}=\frac{-5x}{\dr_y f_2}$. 

We use the previous algorithm to compute the set of values of $\mc{R}_D$. It will give us figure~\ref{res:fig:algo}. 

\smallskip

\etape{First step} 

\begin{itemize}
\item We set $g=f_y'$, which induces a non zero divisor in $\co_D$. Since $\restriction{f_y'}{D_1}=f_2\frac{\dr f_1}{\dr y}$ and $\restriction{f_y'}{D_2}=f_1\frac{\dr f_2}{\dr y}$, one can find an expression of the residues along $D_1$ and $D_2$ with denominator $f_y'$.%, we also have $\rho_{11}=1$ and $\rho_{12}=\frac{f_2(-3x)}{f_y'}$. 
\item We check that $\lrp{\lra{\rho_{11},\rho_{12}},\lra{x,y}}$ is a standard basis of $\mc{R}_{D_1}$ as an $\co_{D_1}$-module. In particular, we have $\val_1(\mc{R}_{D_1})=\lra{-7,-4,-2,-1}\cup\N$. 
\item Similarly, $\lrp{\lra{\rho_{21},\rho_{22}},\lra{x,y}}$ is a standard basis of $\mc{R}_{D_2}$ as an $\co_{D_2}$-module. We then have $\val_2(\mc{R}_{D_2})=\lra{-7,-4,-2,-1}\cup \N$. 
\item From the restrictions of $\rho_1$ and $\rho_2$ to the first branch and the fact that $(x=t_1^3,y=t_1^5)$ is a standard basis of $\co_{D_1}$, we deduce that: 
$\val_1(\mc{R}_D)=\lra{-10,-7,-6,-5,-4,-3,-2,-1}\cup\N.$

In particular, $(\lra{\rho_1,\rho_2},\lra{x,y})$ is a standard basis of the $\co_{D_1}$-module $\restriction{\mc{R}_D}{D_1}$. 
\end{itemize}

%\begin{remar}
%We also need the restriction of $\rho_1$ and $\rho_2$ to the second branch.
%\end{remar}

\etape{Second step}

\s

We set $\mc{M}_0=\lra{(0,-7), (0,-4), (0,-2), (0,-1), (0,0)}$ and $H_0=\lrp{\rho_{21},\rho_{22}}$.
\begin{itemize}
\item Since $-1\in\val_1(\mc{R}_{D_1})$, we have $\mc{M}_1=\mc{M}_0\cup\lra{(-1,v)\ ;\ v\in\val_2(\mc{M}_0)}$ et $H_1=H_0$.
\item Similarly, $\mc{M}_2=\mc{M}_1\cup\lra{(-2,v)\ ;\ v\in\val_2(\mc{M}_1)}$ and $H_2=H_1$.
\item We have $-3\in\val_1(\mc{R}_D)$ but $-3\notin\val_1(\mc{R}_{D_1})$. We have $\val(x\rho_1)=(-3,-5)$. Since $-5\notin\val_2(\mc{M}_2)$, we have: $\mc{M}_3=\mc{M}_{-2}\cup\lra{(-3,-5)}\cup\lra{(-3,-7)}$ et $H_3=H_2\cup\lra{x\rho_1}$. 

By iterating the method, we obtain: 
\item $\mc{M}_{-4}=\mc{M}_{-3}\cup\lra{(-4,v_2)\ ;\ v_2\in\val_2(\mc{M}_{-3})}$ and $H_4=H_3$
\item $\mc{M}_{-5}=\mc{M}_{-4}\cup\lra{(-5,-3)}\cup\lra{(-5,-4)}\cup\lra{(-5,-5)}\cup\lra{(-5,-7)}$ and $H_5=H_4\cup\lra{y\rho_2}$
\item $\mc{M}_{-6}=\mc{M}_{-5}\cup\lra{(-6,-10)}$ and $H_6=H_5\cup\lra{\rho_1}$
\item $\mc{M}_{-7}=\mc{M}_{-6}\cup\lra{(-7,v_2)\ ;\ v_2\in\val_2(\mc{M}_{-6})}$ and $H_7=H_6$
\item $\mc{M}_{-8}=\mc{M}_{-7}$ and $H_8=H_7$
\item $\mc{M}_{-9}=\mc{M}_{-8}$ and $H_9=H_8$
\item $\mc{M}_{-10}=\mc{M}_{-9}\cup\lra{(-10,-6)}\cup\lra{(-10,-7)}\cup\lra{(-10,-10)}$ and $H_{10}=H_9\cup\lra{\rho_2}$. 
\end{itemize}

\begin{figure}[H]
\begin{center}
\begin{tikzpicture}[scale=0.9]
%\shorthandoff{:} 

\draw [step=0.5, gray, very thin] (-5,-5) grid (0.25,0.25);
\draw [->, thick] (-5.5,0) -- (0.6,0) node[at end, below] {{\small $\val_1$}};
\draw [->, thick] (0,-5.5) -- (0,0.5) node[at end, left] {{\small $\val_2$}};
\croix{0}{0}
\croix{0}{-0.5}
\croix{0}{-1}
\croix{0}{-2}
\croix{0}{-3.5}
\croix{-0.5}{0}
\croix{-0.5}{-0.5}
\croix{-0.5}{-1}
\croix{-0.5}{-2}
\croix{-0.5}{-3.5}
\croix{-1}{0}
\croix{-1}{-0.5}
\croix{-1}{-1}
\croix{-1}{-2}
\croix{-1}{-3.5}
\croix{-1.5}{-2.5}
\croix{-1.5}{-3.5}
\croix{-2}{0}
\croix{-2}{-0.5}
\croix{-2}{-1}
\croix{-2}{-2}
\croix{-2}{-3.5}
\croix{-2}{-2.5}
\croix{-2.5}{-1.5}
\croix{-2.5}{-2}
\croix{-2.5}{-2.5}
\croix{-2.5}{-3.5}
\croix{-3}{-5}
\croix{-3.5}{0}
\croix{-3.5}{-0.5}
\croix{-3.5}{-1}
\croix{-3.5}{-2}
\croix{-3.5}{-3.5}
\croix{-3.5}{-2.5}
\croix{-3.5}{-5}
\croix{-3.5}{-1.5}
\croix{-5}{-5}
\croix{-5}{-3}
\croix{-5}{-3.5}

\draw (-5,0.05) node[above]{$-10$};
\draw (-2.5,0.05) node[above]{$-5$};
\draw (0.15,-2.5) node[right]{$-5$};
\draw (0.15,-5) node[right]{$-10$};
\draw (0.25,0.25) node{$0$};
\end{tikzpicture}
\caption{Values of $\mc{R}_D$}
\label{res:fig:algo}
\end{center}
\end{figure}

\end{ex}

\bibliographystyle{alpha-fr}
\bibliography{bibli2}

\end{document}